\documentclass[reqno]{amsart}
\usepackage{epsfig,amsmath,amsfonts,latexsym}
\usepackage{amsthm}
\usepackage[abs]{overpic}
\usepackage[usenames,dvipsnames]{xcolor}
\usepackage{palatino}
\usepackage[hyperfootnotes=false]{hyperref}
\usepackage{caption}
\usepackage{dsfont}
\usepackage{comment}
\usepackage{accents}

\usepackage{mathtools}
\usepackage{cite}
\usepackage{tikz-cd}
\hypersetup{
  colorlinks,
  citecolor=Blue,
  linkcolor=Black,
  urlcolor=Green}
  
\theoremstyle{plain}

\newtheorem{theorem}{Theorem}

\newtheorem{dummy}{anything}[section]
\newtheorem*{thm}{Theorem}
\newtheorem{lemma}[dummy]{Lemma}

\newtheorem{proposition}[dummy]{Proposition}

\newtheorem{corollary}[dummy]{Corollary}

\theoremstyle{definition}
\newtheorem{definition}[dummy]{Definition}

\newtheorem{example}[dummy]{Example}

\newtheorem{remark}[dummy]{Remark}

\theoremstyle{remark}
\newcommand{\N}{\mathbb{N}}
\newcommand{\Z}{\mathbb{Z}}
\newcommand{\R}{\mathbb{R}}
\newcommand{\C}{\mathbb{C}}

\DeclareMathOperator{\crit}{crit}

\DeclareMathOperator{\rk}{rk}
\DeclareMathOperator{\id}{id}

\DeclareMathOperator{\ind}{ind}
\DeclareMathOperator{\atan}{atan}

\def\R{\mathbb{R}}
\def\Z{\mathbb{Z}}
\def\C{\mathbb{C}}
\def\N{\mathbb{N}}

\title{A generalized Poincar\'e--Birkhoff theorem}
\author{Agustin Moreno, Otto van Koert}

\address[A.\ Moreno]{Department of Mathematics \\ Uppsala University \\  Uppsala  \\ Sweden}

\email{agustin.moreno2191@gmail.com}

\address[O.\ van Koert]{Department of Mathematical Sciences \\ Seoul National University \\  Seoul  \\ South Korea}

\email{okoert@snu.ac.kr}

\date{}

\begin{document}

\maketitle

{\centering\footnotesize  \textit{To H.\ Poincar\'e, who taught us much;\\ To A. Floer, who followed suit;\\ To C.\ Viterbo, now on his 60th birthday, who took the cue; \\ and to all those who stand on the Shoulders of Giants.}\par}

\begin{abstract}
We prove a generalization of the classical Poincar\'e--Birkhoff theorem for Liouville domains, in arbitrary even dimensions. This is inspired by the existence of global hypersurfaces of section for the spatial case of the restricted three-body problem \cite{MvK}.
\end{abstract}

\tableofcontents

\section{Introduction} 

\textbf{Poincar\'e--Birkhoff theorem, and the planar restricted three-body problem.} The problem of finding closed orbits in the planar case of the restricted three-body problem goes back to ground-breaking work in celestial mechanics of Poincar\'e \cite{P87,P12}, building on work of G.W.\ Hill on the lunar problem \cite{H78}. The basic scheme for his approach may be reduced to:
\begin{itemize}
    \item[(1)] Finding a \emph{global surface of section} for the dynamics;
    \item[(2)] Proving a fixed point theorem for the resulting first return map.
\end{itemize}
This is the setting for the celebrated Poincar\'e--Birkhoff theorem, proposed and confirmed in special cases by Poincar\'e and later proved in full generality by Birkhoff in \cite{Bi13}. The statement can be summarized as: if $\tau: A\rightarrow A$ is an area-preserving homeomorphism of the annulus $A=[-1,1]\times S^1$ that satisfies a \emph{twist} condition at the boundary, then it admits infinitely many periodic points of arbitrary large period. 

In \cite{MvK}, the authors proved the existence of $S^1$-families of global hypersurfaces of section for the \emph{spatial} restricted three-body problem (in the low-energy range, i.e.~below and slightly above the first critical value, and independent of mass ratio), fully and non-perturbatively generalizing step (1) in the above approach to the spatial situation. 
The relevant return map $\tau$ is a Hamiltonian diffeomorphism defined on the interior of the global hypersurface of section, which is symplectomorphic to the interior of a Liouville domain $(\mathbb{D}^*S^2,\omega)$, where $\omega$ is deformation equivalent to the standard symplectic form.
Furthermore $\tau$ extends smoothly to the boundary of the global hypersurface of section, and gives rise to a homeomorphism of $(\mathbb{D}^*S^2,\omega)$ that is a Hamiltonian diffeomorphism on the interior.
Drawing inspiration from this situation, in this paper, we propose a general fixed-point theorem for Liouville domains, as an attempt to address step (2) for the spatial case.

\subsection*{Fixed-point theory of Hamiltonian twist maps} The periodic points of $\tau$ are either boundary periodic points, which give planar orbits, or interior periodic points which are in 1:1 correspondence with spatial orbits. We are interested in finding \emph{interior} periodic points.

\medskip

\textbf{The Hamiltonian twist condition.} 
We propose a generalization of the twist condition introduced by Poincar\'e, for the Hamiltonian case and for arbitrary Liouville domains. Let $(W,\omega=d\lambda)$ be a $2n$-dimensional Liouville domain, and consider a Hamiltonian symplectomorphism $\tau$ of $W$. Let $(B,\xi)=(\partial W,\ker \alpha)$ be the contact manifold at the boundary where $\alpha=\lambda\vert_B$, and $R_\alpha$ the Reeb vector field of $\alpha$ (uniquely determined via the equations $d\alpha(R_\alpha,\cdot)=0,\alpha(R_\alpha)=1$). Recall that $\tau$ is Hamiltonian if $\tau=\phi_H^1$, where $\phi_H^t$ is the isotopy of $W$ defined by $\phi_H^0=id$, $\frac{d}{dt}\phi_H^t=X_{H_t}\circ \phi_H^t$, where we write $H_t=H(t,\cdot)$, and $X_{H_t}$ is the Hamiltonian vector field of $H_t$ defined via $i_{X_{H_t}}\omega=-dH_t$. The Liouville vector field $V_\lambda$ is defined via $i_{V_\lambda}\omega=\lambda$. 

\begin{definition}(\textbf{Hamiltonian twist map})\label{def:twistmap}
We say that $\tau$ is a \emph{Hamiltonian twist map} (with respect to $\alpha$), if $\tau$ is generated by a \emph{smooth} Hamiltonian $H:W \times \mathbb{R} \rightarrow \mathbb{R}$ which satisfies $X_{H_t}\vert_B=h_tR_\alpha$ for some \emph{positive} and smooth function $h:B \times \mathbb{R} \rightarrow \mathbb{R}^+$.
\end{definition}

\begin{remark}
For the purposes of this article, one may relax the smoothness assumption on $H$ to $C^2$ regularity.
\end{remark}

In particular, $H_t\vert_B \equiv const$ on $B$, and $\tau(B)\subset B$. We have $h_t=dH_t(V_\lambda)\vert_B$ is the derivative of $H_t$ in the Liouville direction $V_\lambda$ along $B$, which we assume strictly positive. Also, $\tau\vert_B$ is the time-$1$ map of a positive reparametrization of the Reeb flow on $B$. But note that, while the latter condition is only localized at $B$, the twist condition is of a \emph{global} nature, as it requires global smoothness of the generating Hamiltonian (cf.\ \cite[Rk.\ 1.4]{MvK}). 

\medskip

Here is a simple example illustrating why the smoothness of the Hamiltonian is relevant for the purposes of fixed points:

\begin{example}[Integrable twist maps]\label{ex:inttwistmap} Consider $M=S^n$, $n\geq 1$ with its round metric and its cotangent bundle $T^*M=\{(q,p)\in \mathbb{R}^{2n+2}:\langle q,p\rangle=0,\; \vert q \vert =1\}$. Let $H: T^*M \rightarrow \mathbb{R}$, $H(q,p)=2\pi|p|$ (\emph{not} smooth at the zero section); $\phi_H^1$ extends to all of $\mathbb{D}^*M$ as the identity. It is a positive reparametrization of the Reeb flow at $S^*M$, generating a full turn of the geodesic flow, and all orbits are fixed points with fixed period. If we smoothen $H$ near $|p|=0$ to $K(q,p)=2\pi g(|p|)$, with $g(0)=g^\prime(0)=0$, then $\tau=\phi_K^1: \mathbb{D}^*M \rightarrow \mathbb{D}^*M$, $\tau(q,p)=\phi_H^{2\pi g^\prime(|p|)}(q,p)$, is now a Hamiltonian twist map. If $g^\prime(|p|)=l/k\in \mathbb{Q}$ with $l,k$ coprime, then $\tau$ has a simple $k$-periodic orbit; therefore $\tau$ has simple interior orbits of arbitrary large period (cf.\ \cite[p.\ 350]{KH95}, \cite{M86}, for the case $M=S^1$).
\end{example}

The Hamiltonian twist condition will be used to extend the Hamiltonian to a Hamiltonian that is admissible for computing symplectic homology. The extended Hamiltonian can have additional $1$-periodic orbits and these, as well as $1$-periodic orbits on the boundary, need be distinguished from the interior periodic points of $\tau$. We impose the following conditions to do so.

\medskip

\textbf{Index growth.} We consider a suitable index growth condition on the dynamics on the boundary, which is satisfied in the restricted three-body problem whenever the \emph{planar} dynamics is strictly convex (see Thm.\ \ref{lemma:indexposapp}). This assumption will allow us to separate boundary and extension orbits from interior ones via the index.

\medskip

We call a strict contact manifold $(Y,\xi=\ker \alpha)$ \emph{strongly index-definite} if the contact structure $(\xi,d\alpha)$ admits a symplectic trivialization $\epsilon$ with the property that
\begin{itemize}
    \item There are constants $c>0$ and $d\in \mathbb{R}$ such that for every Reeb arc\footnote{We will refer to the restriction of a Reeb orbit or Hamiltonian orbit to a finite interval as a Reeb arc or Hamiltonian arc.} $\gamma:[0,T]\rightarrow Y$ of Reeb action $T=\int_0^T \gamma^*\alpha$ we have
    $$
    \vert\mu_{RS}(\gamma;\epsilon)\vert\geq c T+d,
    $$
    where $\mu_{RS}$ is the Robbin--Salamon index \cite{RS93}.
\end{itemize}

Index-positivity is defined similarly, where we drop the absolute value. A variation of this notion was explored in Ustilovsky's thesis \cite{U99}. 
He imposed the additional condition $\pi_1(Y) = 0$.
With this extra assumption, the concept of index positivity becomes independent of the choice of trivialization, although the exact constants $c$ and $d$ still depend on the trivialization $\epsilon$. The global trivialization will be important when considering extensions of our Hamiltonians, as it will allow us to measure the index growth of potential new orbits.

\medskip

\textbf{Fixed-point theorems.} We propose the following generalization of the Poincar\'e--Birkhoff theorem:
\begin{theorem}[Generalized Poincar\'e--Birkhoff theorem]
\label{thmA}
Suppose that $\tau$ is an exact symplectomorphism of a connected Liouville domain $(W,\lambda)$, and let $\alpha=\lambda\vert_B$. Assume the following:
\begin{itemize}
\item\textbf{(Hamiltonian twist map)} $\tau$ is a Hamiltonian twist map, where the generating Hamiltonian is at least $C^2$. In addition, assume all fixed points of $\tau$ are isolated;
\item\textbf{(Index-definiteness)} If $\dim W\geq 4$, then assume $c_1(W)\vert_{\pi_2(W)}=0$, and $(\partial W, \alpha)$ is strongly index-definite; 
\item\textbf{(Symplectic homology)} $SH_\bullet(W)$ is infinite dimensional.
\end{itemize}
Then $\tau$ has simple interior periodic points of arbitrarily large (integer) period.
\end{theorem}

\begin{remark}\label{rk:thmA} Let us discuss some aspects of the theorem:

\begin{enumerate}

    \item\textbf{(Grading)} We impose the assumptions $c_1(W)\vert_{\pi_2(W)}=0$ (i.e.\ $W$ is symplectic Calabi-Yau) to have a well-defined integer grading on symplectic homology.
    
    \smallskip
    
    \item\textbf{(Surfaces)} If $\dim W=2$, then the condition that $SH_\bullet(W)$ is infinite dimensional just means that $W$ is not $D^2$ (see App.\ \ref{app:SH_surfaces}); for $D^2$ we have $SH_\bullet(D^2)=0$, and a rotation on $D^2$ gives an obvious counterexample to the conclusion. 
    In the surface case, the argument simplifies, and one can simply work with homotopy classes of loops rather than the grading on symplectic homology. The Hamiltonian twist condition implies the classical twist condition for $W=\mathbb{D}^*S^1$, due to orientations.
  
   \smallskip
    
    \item\textbf{(Cotangent bundles)} The symplectic homology of the cotangent bundle of a closed manifold with finite fundamental group is well-known to be infinite dimensional, due to a result of Viterbo \cite{V99, V18} (see also \cite{AS,SW06}), combined e.g.\ with a theorem of Gromov \cite[Sec.\ 1.4]{G78}. We have $c_1(T^*M)=0$ whenever $M$ is orientable. 
As for the existence of a global trivialization of the contact structure $(\xi,d\lambda_{can})$, we note the following:
\begin{itemize}
\item if $\Sigma$ is an oriented surface, then $S^*\Sigma$ admits such a global symplectic trivialization;
\item if $M^3$ is an orientable $3$-manifold, then $S^* M^3$ also admits such a global symplectic trivialization;
\item  In addition, we know that symplectic trivializations of the contact structure on $(S^*S^2,\lambda_{can})$ are unique up to homotopy, since $[S^*S^2,Sp(2)]\cong H^1(S^*S^2;\Z)=0$.
\end{itemize}
    
     \smallskip
    
    \item\textbf{(Fixed points)} If fixed points are non-isolated, then we vacuously obtain infinitely many of them, although we cannot conclude that their periods are unbounded; ``generically'', one expects finitely many fixed points.
    
     \smallskip
    
    \item\textbf{(Long orbits)} If $W$ is a global hypersurface of section for some Reeb dynamics, with return map $\tau$, interior periodic points with long (integer) period for $\tau$ translate into spatial Reeb orbits with long (real) period; see Lemma~\ref{lemma:appendix}.
    
     \smallskip
    
    \item\textbf{(Katok examples)} There are well-known examples due to Katok \cite{K73} of Finsler metrics on spheres with only finitely many simple geodesics, which are arbitrarily close to the round metric (we review them in App.\ \ref{app:Katok}); they admit global hypersurfaces of section with Hamiltonian return maps, for which the index-definiteness and the condition on symplectic homology hold. It follows that the return map does not satisfy the twist condition for any choice of Hamiltonians.
    
     \smallskip
    
    \item\textbf{(Spatial restricted three-body problem)} From the above discussion and \cite{MvK}, we gather: the only standing obstruction for applying the above result to the spatial restricted three-body problem, in case where the planar problem is strictly convex, is the Hamiltonian twist condition. Here, note that symplectic homology is invariant under deformations of Liouville domains; see e.g.\ \cite{BR} for a paper with detailed proofs. 
    This would give a proof of existence of \emph{spatial} long orbits in the spirit of Conley \cite{C63}, which could in principle be collision orbits. Since the geodesic flow on $S^2$ arises as a limit case (i.e.\ the Kepler problem), it should be clear from the discussion on Katok examples that this is a subtle condition. In \cite{MvK}, we have computed a generating Hamiltonian for the integrable case of the rotating Kepler problem; it does \emph{not} satisfy the twist condition in the spatial case (in the planar case, a Hamiltonian twist map was essentially found by Poincar\'e). This does not mean a priori that there is not \emph{another} generating Hamiltonian which does, but this seems rather unlikely.
\end{enumerate}

\end{remark}

As a particular case of Thm.\ \ref{thmA}, we state the above result for star-shaped domains in cotangent bundles, as of independent interest (cf.\ \cite{H11}): 

\begin{theorem}
\label{thmB}
Suppose that $W$ is a fiber-wise star-shaped domain in the Liouville manifold $(T^*M,\lambda_{can})$, where $M$ is simply connected, orientable and closed, and assume that $\tau:W\to W$ is a Hamiltonian twist map. If the Reeb flow on $\partial W$ is strongly index-positive, and if all fixed points of $\tau$ are isolated, then $\tau$ has simple interior periodic points of arbitrarily large period. 
\end{theorem}

The above holds in particular for $M=S^2$, as explained in Remark \ref{rk:thmA} (3). One difference with \cite{H11} is that we work with compact domains in cotangent bundles and conclude that periodic points are interior, at the expense of imposing index-positivity.

\medskip

\textbf{Sketch of the proof.} The proof is fairly simple: due to the twist condition we can extend the map $\tau$ to a Hamiltonian diffeomorphism $\widehat \tau$ that is generated by a weakly admissible Hamiltonian (defined in Section~\ref{sec:proof_gpb}). This allows us to appeal to symplectic homology.
In particular, we will show $\varinjlim_k HF_\bullet(\widehat \tau^k) = SH_\bullet(W)$ (Lemma \ref{lemma:extension__Hamiltonian_linear}). Using an index filtration (via index-definiteness \emph{and} the twist condition), we can show that all generators contributing to homology are actually fixed points of some $\tau^k$, rather than fixed points of the extension. The crucial technical input is Lemma~\ref{indexpositivitylemma}. If the minimal periods of periodic points of $\tau$ are bounded, then we can show using a spectral sequence, involving local Floer homology groups, that the rank of the resulting symplectic homology should also be bounded, leading to a contradiction. Alternatively, one could use the methods used for the proof of the Conley conjecture \cite{G10,H11} to finish the proof.

\medskip

\textbf{Remarks on the twist condition and generalizations.} If the Liouville domain is a surface, this definition of the Hamiltonian twist condition is not restrictive, and implements the idea sketched above in a simple way.
In higher dimensions, the Hamiltonian twist condition is much more restrictive. Some examples illustrating the nature of the twist condition and applications of the above theorem will be presented in Appendix~\ref{sec:examples}. Given the above sketch of the proof, there is obviously some freedom in Def.\ \ref{def:twistmap} that allows the same methods to work. For example, if the vector field $X_{H_t}$ is sufficiently $C^1$-close to a positive reparametrization of the Reeb vector field, then the methods will still go through. However, we will not pursue this generalization because its depends on details that make the formulation awkward and difficult to check.
We list some other generalizations, whose proofs will not be worked out in detail:
\begin{itemize}
\item (\textbf{Action positivity}) One can impose constraints on the functions  $h_t$ in the Hamiltonian twist condition that force the periodic orbits in the extension to have large action under iterates.
In the setting of cotangent bundles, one can then use a theorem of Gromov \cite[Sec.\ 1.4]{G78} cited below, to construct infinitely many interior periodic points.

\smallskip

\item (\textbf{Isolated sets}) The assumption that the fixed points are isolated can be replaced by the weaker assumption that the fixed point set consists of a finite union of submanifolds. This is based on a slight generalization of local Floer homology, and is useful when studying integrable systems and their perturbations.

\smallskip

\item (\textbf{Non-vanishing symplectic homology}) The condition $\dim SH_\bullet(W) =\infty$ can be replaced by the condition $SH_\bullet(W) \neq 0$.
The key point here is that non-vanishing symplectic homology implies its unit is non-trivial.
Then the methods of the proof of the Conley conjecture \cite{G10, H11} can be applied to conclude the existence of infinitely many simple periodic points. Strong index-definiteness is needed to show that these periodic points do not correspond to boundary and extension orbits, and so are interior.
\end{itemize}
\begin{remark}
Concerning the last generalization, we remark that we don't know a single example of a Liouville domain $(W,\lambda)$ with $c_1(W)=0$, $SH_\bullet(W) \neq 0$, and $\dim SH_\bullet(W)<\infty$. 
\end{remark}

\medskip

\textbf{Acknowledgements.} The authors thank Urs Frauenfelder, for suggesting this problem to the first author, for his generosity with his ideas and for insightful conversations throughout the project; Lei Zhao, Murat Saglam, Alberto Abbondandolo, and Richard Siefring, for further helpful inputs, interest in the project, and discussions. We are also grateful to an anonymous referee for careful comments on earlier versions. The first author has also significantly benefited from several conversations with Kai Cieliebak in Germany and Sweden, as well as with Alejandro Passeggi in Montevideo, Uruguay. This research was started while the first author was affiliated to Augsburg Universit\"at, Germany. The first author is also indebted to a Research Fellowship funded by the Mittag-Leffler Institute in Djursholm, Sweden, where this manuscript was finalized.
The second author was supported by NRF grant NRF-2019R1A2C4070302, which was funded by the Korean Government.

\section{Motivation and background}
\label{sec:prelim}

\subsection*{Hypersurfaces of section, return maps, and open books}

\begin{definition}
Suppose that $Y$ is a compact, oriented, smooth manifold with a non-singular autonomous flow $\phi_t$.
We call an oriented, compact hypersurface $\Sigma$ in $Y$ a \emph{global hypersurface of section} for $\phi_t$ if
\begin{itemize}
\item the set $\partial \Sigma$ is an invariant set for the flow $\phi_t$ (if non-empty);
\item the flow $\phi_t$ is positively transverse to the interior of $\Sigma$;
\item for all $x\in Y \setminus \partial \Sigma$ there are $t_+>0$ and $t_-<0$ such that $\phi_{t_+}(x) \in \Sigma$ and $\phi_{t-}(x)\in \Sigma$.
\end{itemize} 
\end{definition}
Given a global hypersurface of section we can define a return map $\tau$ as follows: for each $x \in \mbox{int}(\Sigma)$ we choose a minimal $t_{+}(x)>0$ as in the definition above. Then we put $\tau(x) = \phi_{t_+(x)}(x)$. Periodic points of $\tau$ then correspond to closed orbits of $\phi_t$. In general, there is no continuous extension to the boundary, although it is unique whenever exists. Although global hypersurfaces of section do not have good stability properties in higher dimensions, we found that they can be constructed in certain classes of Hamiltonian dynamical systems that admit an involution. This class includes the restricted three-body problem and several variations (e.g.\ suitable Stark-Zeeman systems \cite{MvK}). 

This notion is also closely related to the notion of an open book decomposition. This consists of a fiber bundle $\pi: Y \backslash B \rightarrow S^1$, where $B\subset Y$ is a codimension $2$ submanifold with trivial normal bundle (called the \emph{binding}), such that $\pi$ coincides with the angular coordinate along some choice of collar neighbourhood $B\times \mathbb{D}^2$ of $B$. The \emph{pages} of the open book are the closure of the fibers of $\pi$, all having $B$ as boundary. Whenever $\phi_t$ is a Reeb dynamics of a contact form $\alpha$ on $Y$ which is adapted to the open book (i.e.\ $\alpha\vert_B$ is also contact, and $d\alpha$ is symplectic on the pages), each page is a global hypersurface of section, and the return map preserves the symplectic form $d\alpha$. This is precisely the situation in \cite{MvK}.

In App.\ \ref{app:longorbits}, we will collect some standard facts which apply for return maps arising from Reeb dynamics, as described here, for which Thm.\ A may be applied.

\section{Preliminaries on symplectic homology}
\subsection{Liouville domains and Hamiltonian dynamics}
There are various forms of Hamiltonian Floer homology for Liouville domains: these are all referred to as \emph{symplectic homology}.
The first version was due to Floer-Hofer, \cite{FH}. See also section~5 of \cite{H89} for an even earlier version, called symplectology.
However, we will review the version due to Viterbo, \cite{V99, V18}. Roughly speaking, this is a ring with unit that encodes both topological and dynamical data; it is the homology of a chain complex that is freely generated by $1$-periodic Hamiltonian orbits.

We now fix conventions. Consider a \emph{Liouville domain} $(W,\lambda)$, i.e.\ $(W,d\lambda)$ is a compact symplectic manifold with boundary, and the vector field $X$ defined by the equation $\iota_X d\lambda =\lambda$ is outward pointing along each boundary component of $W$.
This vector field is the \emph{Liouville vector field}.
The $1$-form $\lambda$ is the \emph{Liouville form}, and its restriction to $\partial W$, which we denote by $\alpha$, is a contact form.

Given a Liouville domain $(W,\lambda)$ we build its completion to a \emph{Liouville manifold} by attaching a cylindrical end: 
$$
(\widehat W, \widehat \lambda):= (W,\lambda) \cup_\partial ([1,\infty) \times \partial W,r \alpha ).
$$
Throughout the paper we will consider smooth functions of the form $H:W\times S^1 \to \R$, a (time-dependent) \emph{Hamiltonian} on $W$. 
Given such a Hamiltonian, we define its Hamiltonian vector field $X_H$ via
$$
\iota_{X_H} d\lambda = -dH.
$$
We denote the set of $1$-periodic orbits of $X_H$ by $\mathcal P(H)$.
For the purpose of Floer theory on non-compact manifolds we will need a suitable class of Hamiltonians to work with. First, we recall the spectrum of a contact form $\alpha$.
If $\mathcal P(\alpha)$ denotes the set of all periodic Reeb orbits (including covers and without period bound), then 
$$
\mbox{spec}(\alpha)
=\{ a \in \R ~|~\text{there is }\gamma \in \mathcal P(\alpha) \text{ such that } a= \mathcal A(\gamma) \},
$$
where the \emph{action} is defined as $\mathcal{A}(\gamma)=\int_\gamma \alpha$. 
\begin{definition} We recall some standard terminology.
\begin{itemize}
    \item A $1$-periodic orbit $\gamma\in \mathcal P(H)$ is \emph{non-degenerate} if $dFl^{X_H}_1(\gamma(0))-\id$ invertible.\\
    \item The Hamiltonian $H$ is \emph{non-degenerate} if all $\gamma\in \mathcal P(H)$ are non-degenerate.\\
    \item A Hamiltonian $ H$ on $\widehat W$ is \emph{linear at infinity} if at the cylindrical end $ H$ has the form $ H(r,b,t) =cr+d$ for some constants $c>0$ and $d$. In this case we write slope$( H):=c$.
    \item A Hamiltonian $ H$ that is non-degenerate and linear at infinity with slope$(H) \notin \mbox{spec}(\alpha)$ will be called \emph{admissible}.
\end{itemize}
We call an $S^1$-family of almost complex structures $J=J_t$ on a Liouville manifold $(\widehat W,\widehat \lambda)$ \emph{SFT-like} if
\begin{itemize}
\item it is compatible with $(T\widehat W,d\widehat \lambda)$; and
\item on the cylindrical end it satisfies $\mathcal L_{r\partial r} J=0$, $J \xi =\xi$, and $Jr\partial_r =R_\alpha$.
\end{itemize}
We denote by $\mathcal{J}$ the space of such families of almost complex structures.
\end{definition}
\subsection{Conley-Zehnder index, Robbin-Salamon index and mean index}
\label{sec:def_indices}
We will also need invariants of Hamiltonian orbits, i.e.\ the Conley-Zehnder index, or more generally, the Robbin-Salamon index, and the mean index. Assume that $x:\mathbb{R} \to \widehat W$ is an orbit of $X_H$. Take a symplectic trivialization $\epsilon:\mathbb{R} \times \R^{2n} \to x^*T\widehat W,~(t,v)\mapsto \epsilon_t(v) \in T_{x(t)} \widehat W$. 
Then we get a path of symplectic matrices associated with $x$, namely $\psi_t = \epsilon_t^{-1} \circ dFl^{X_H}_t \circ \epsilon_0.$ We can then define the Robbin-Salamon index of $x$ as $\mu_{RS}(x\vert_{
[0,T]},\epsilon) :=
\mu_{RS}(\psi|_{[0,T]})$. If $\psi_T-\id$ is invertible, then the Robbin-Salamon index reduces to the Conley-Zehnder index.
The case of Reeb flows is done similarly; we simply restrict the linearized Reeb flow to the symplectic vector bundle $(\xi,d\alpha)$. Similarly, we define the mean index of a $1$-periodic orbit $x$ as $
\Delta(x,\epsilon):=\Delta(\psi)$, 
where $\Delta(\psi)$ is the mean index of the symplectic path $\psi$.\footnote{A description of the mean index can be found on page 1318 of of \cite{SZ}.}
We have the following properties (see e.g.\ section~ 3.1.1 of  \cite{GG15}):
$$
\begin{array}{cc}
   (1)  & \vert \mu_{RS}(x\vert_{[0,T]},\epsilon)-\Delta(x,\epsilon)\vert \leq \frac{\dim W}{2}, \mbox{ for all }T; \\
   (2)  & \lim_{T\rightarrow +\infty}\frac{\mu_{RS}(\psi|_{[0,T]},\epsilon)}{T}=\Delta(x,\epsilon);\\
  (3) & \Delta(x^{(k)},\epsilon)=k\Delta(x,\epsilon),
\end{array}
$$
where we interpret the $k$-fold catenation $x^{(k)}$, a $k$-periodic orbit of $H$, as a $1$-periodic orbit of the iterated Hamiltonian $H^{\# k}$.

\begin{definition}
We will call a Hamiltonian flow on $W$ \emph{strongly index-definite} if there is a symplectic trivialization $\epsilon_W:W\times \R^{2n} \to TW$, and constants $c>0, d$ and such that for every orbit of $X_H$, we have
$$
|\mu_{RS}(x\vert_{[0,T]},\epsilon)|\geq c T+d.
$$
\end{definition}
The notion of \emph{strong index-positivity} is obtained by dropping the absolute value in the above definition, and similarly for \emph{strong index-negativity}. As in the Introduction, we can also define it for Reeb flows. Here are some examples:
\begin{lemma}
Suppose that $(M,g)$ is a closed Riemannian manifold with positive sectional curvature.
Assume in addition that the contact structure $(S^*M,(\xi,d\alpha))$ admits a global symplectic trivialization.
Then $(S^*M,d \alpha)$ is strongly index-positive. 
\end{lemma}
Other examples are complements of Donaldson hypersurfaces in monotone symplectic manifolds provided the degree is sufficiently high and symplectically trivial: these manifolds are index negative.

\subsection{Hamiltonian Floer homology and symplectic homology}
Given \emph{Floer data} $(J, H)$ of an SFT-like $J$ and an admissible $H$, we note the following:
\begin{itemize}
\item There are no $1$-periodic orbits of $X_{ H}$ on the cylindrical end, because of the spectrum assumption.
\item Non-degenerate $1$-periodic orbits of $X_{ H}$ are isolated.
\end{itemize}
Then $\mathcal P( H)$ consists of finitely many $1$-periodic orbits. Informally speaking, we think of Floer homology as ``Morse homology'' of the following action functional:
\[
\mathcal A_{ H}: W^{1,2}(S^1=\R /\Z, \widehat W) \longrightarrow \R,
\quad
\gamma \longmapsto \int_{S^1} \gamma^* \widehat \lambda -\int_0^1  H(\gamma(t),t)dt.
\]
This functional has the property $\mathcal A_{H^{\# k}}(x^{(k)}) = k \mathcal A_{H}(x)$ for iterates. A computation shows that $\crit {\mathcal A}_{ H} =\mathcal P( H)$, and we define the Floer chain complex as:
$$
CF_{\bullet}(\widehat W, \widehat \lambda, H, J)
:=\bigoplus_{\gamma \in \mathcal P( H)} \Z_2 \langle \gamma \rangle.
$$
We grade this chain complex by the Conley-Zehnder index, so $\deg \gamma:=\mu_{CZ}(\gamma,\epsilon)$.
We make a couple of comments:
\begin{itemize}
\item in the standard procedure, we choose a capping disk $\tilde \gamma$ of a contractible $1$-periodic orbit $\gamma$, and a symplectic trivialization $\epsilon_{\tilde \gamma}$ of $\tilde \gamma^*T\widehat W$. This gives a trivialization $e_t: (\R^{2n},\omega_0) \to T_{\gamma(t)}\widehat W$ by restriction.
We then define the Conley-Zehnder index of an orbit as in Section~\ref{sec:def_indices}.
Once the capping disk $\tilde \gamma$ is fixed, this index is independent of the choice of trivialization on a \emph{fixed} capping disk. The index does depend on the choice of capping disk but not if $c_1(W)\vert_{\pi_2(W)}=0$.
\item for non-contractible orbits, one needs to choose a reference loop $c$ and a reference symplectic trivialization $\epsilon_c$ for each free homotopy $[c]\in \tilde \pi_1(W)$.
Given a $1$-periodic orbit $x$ in the same free homotopy class as $c$ we choose a connecting cylinder $S$; the trivialization extends over $S$, and we can then define the Conley-Zehnder index as before. 
\item we choose the simpler, but more restrictive approach to use a global symplectic trivialization on some subdomain $\tilde W$. The existence of such a trivialization implies that $c_1(T\tilde W)=0$.
This approach obviously reduces to the previous approach provided that capping disks or connecting cylinders can be chosen to lie in the domain of definition of the global symplectic trivialization.
\end{itemize}
If we define an $L^2$-metric on $W^{1,2}(S^1,x^*T\widehat W)$ by
$$
\langle X, Y \rangle =
\int_{0}^1 \omega(X(t), J_t(x(t))Y(t) \, )dt,
$$
then the Floer equation is the $L^2$-gradient ``flow''\footnote{The flow is strictly speaking not defined, since it leads to an ill-posed initial value problem.} of the above functional: for a cylinder $u:Z=\R \times S^1 \to \widehat W$, this is
\begin{equation}
\label{eq:Floer}
(du-X_{ H}\otimes dt)^{0,1}=0,\;
\lim_{s\to \pm \infty}u(s,t) =x_{\pm}(t).
\end{equation}
Solutions to this equation are called \emph{Floer trajectories.} Given $1$-periodic orbits $x_+, x_- \in \mathcal P( H)$, the moduli space of Floer trajectories is
\[
\mathcal M(x_+,x_-)
:=
\{ u: Z \to \widehat W~|~u\text{ satisfies }\eqref{eq:Floer}
\}.
\]
In general, this space does not need to have a manifold structure. To obtain this extra structure, we first interpret Equation~\eqref{eq:Floer} as a section of a vector bundle, via
\[
\bar \partial_F: \mathcal P(x_+,x_-) \longrightarrow \mathcal E(x_+,x_-),
\quad
u \longmapsto (du-X_{ H} \otimes dt)^{0,1} \in L^p(Z,\Omega^{0,1}(u^*T \widehat W) \,).
\]
Here $\mathcal P(x_+,x_-)$ is a Banach manifold of cylinders of class $W^{1,p}$ that are $W^{1,p}$-pushoffs of smooth cylinders that exponentially converge to the given asymptotes $x_+$ and $x_-$, and $\mathcal E(x_+,x_-)$ is a Banach bundle over $\mathcal P(x_+,x_-)$ whose fiber over $u\in \mathcal P(x_+,x_-)$ is $ L^p(Z,\Omega^{0,1}(u^*T \widehat W) \,)$. For details, see Ch.\ 8 in \cite{AD}.
We will denote the linearization of $\bar \partial_F$ at $u\in \mathcal P(x_+,x_-)$ by $D_u\bar \partial_F$.

\begin{proposition}
For Floer data $(J, H)$ and $u\in \mathcal M(x_+,x_-)$, $D_u\bar \partial_F$ is a Fredholm operator of index 
$$
\ind D_u\bar \partial_F =
\mu_{CZ}(x_+,\epsilon)
-
\mu_{CZ}(x_-,\epsilon),
$$
where $\epsilon$ is a symplectic trivialization of $u^*T\widehat W$.
\end{proposition}
In addition, we can always choose suitable Floer data close to initial Floer data such all moduli space are transversely cut out: 
\begin{proposition}
\label{prop:regular}
There is a dense set $\mathcal J_{reg} \subset \mathcal J$ with the property for all $J\in \mathcal J_{reg}$, the linearized operator $D_u\bar \partial_F$ is surjective for all $u \in \mathcal M(x_+,x_-)$, and so $\mathcal M(x_+,x_-)$ is a smooth manifold of dimension $\mu_{CZ}(x_+,\epsilon)
-
\mu_{CZ}(x_-,\epsilon)$.
\end{proposition}
Floer data $(J, H)$ as in Proposition~\ref{prop:regular} will be called \emph{regular Floer data}. We now have all the basic ingredients in place: choose regular Floer data $(J, H)$, and define the boundary operator for the chain complex $CF_{\bullet}(\widehat W, \widehat \lambda,  H, J)$ via
$$
\partial x_+
=
\sum_{x_- \in \mathcal P( H), \;\deg(x_-)=\deg(x_+)-1}
\#_{\Z_2} \left(\mathcal M(x_+,x_-)/\R \right)\, \cdot x_-.
$$
Here we have modded out $\mathcal M(x_+,x_-)$ by the reparametrization action in the domain, and the resulting quotient spaces can be compactified, so the coefficients in the above sum are actually finite.

\begin{lemma}
This linear map is a differential: $\partial \circ \partial =0$.
\end{lemma}
The Floer homology of $( \widehat W,\widehat \lambda, J, H)$ is then defined as the homology
\[
HF_\bullet(\widehat W,\widehat \lambda,J, H):=
H_\bullet(
CF_{\bullet}(\widehat W,\widehat \lambda,J,H),\partial
)
.
\]
\begin{remark}
In the case of closed symplectic manifolds, Floer homology is independent of the choice of Floer data.
This is not the case for Liouville domains, and this is the next topic we will deal with.
\end{remark}

\subsection{Continuation maps and symplectic homology}
Assume that $H_1$ and $H_2$ are admissible Hamiltonians on a Liouville manifold $\widehat W$ with $\mathrm{slope}(H_1)\leq \mathrm{slope}(H_2)$.
We interpolate between them via
$$
K:\widehat W \times S^1 \times \R \longrightarrow \R,\quad
(w,t,s)\longmapsto K_s(w,t),
$$
where we impose the monotonicity condition $\partial_s K\leq 0$, \footnote{This is the monotonicity condition for the continuation map with our conventions for the Floer equation.
Our conventions agree with those in \cite{F}.} and 
\[
K_s(w,t)=
\begin{cases}
H_1(w,t), & \text{if } s\gg 0 \\
H_2(w,t), & \text{if } s\ll 0.
\end{cases}
\]
We then consider the parametrized Floer equation for $u:Z \to \widehat W$:
\[
(du-X_{K}\otimes dt)^{0,1} =0,\quad
\lim_{s \to \infty} u(s,t) = x_+(t) \in \mathcal P(H_1),
\quad
\lim_{s \to -\infty} u(s,t) = x_-(t) \in \mathcal P(H_2).
\]
The results of the Fredholm theory mentioned in the previous section also apply in this setup, and we can define a \emph{continuation map} as
\[
\begin{split}
c_{12}: CF_\bullet(\widehat W,\widehat \lambda, J, H_1) &\longrightarrow CF_\bullet(\widehat W,\widehat \lambda, J, H_2),
\\
x_+ & \longmapsto \sum_{x_- \in \mathcal P(H_2),\; \deg(x_-)=\deg(x_+)}
\#_{\Z_2} \mathcal M(x_+,x_-, J, K) \cdot x_-,
\end{split}
\]
where $\mathcal M(x_+,x_-, J, K)$ is the moduli space of Floer trajectories of the parametrized Hamiltonian $K$.
\begin{lemma}
The map $c_{12}$ is a chain map, and the induced map on homology is independent of $J,K$.
\end{lemma}
We also write $c_{12}$ for the induced map on Floer homology:
$$
c_{12}: HF_\bullet(\widehat W,\widehat \lambda, J, H_1) \longrightarrow HF_\bullet(\widehat W,\widehat \lambda, J, H_2).
$$
Symplectic homology is then defined as the direct limit over a direct system $\{ H_i \}_i$ of admissible Hamiltonians for whose slopes $\mathrm{slope}(H_i)$ increase to $\infty$,
\begin{equation}
\label{eq:def_SH}
SH_\bullet(W,\lambda, J,\{ H_i \}_i):=
\varinjlim_{c_{ij},~j>i} HF_\bullet(\widehat W,\widehat \lambda, J,H_i).
\end{equation}

\begin{remark}
Symplectic homology is independent of $J$, and the sequence of Hamiltonians $\{ H_i \}_i$.
We will henceforth write $SH_\bullet(W,\lambda)$, or $SH_\bullet(W)$ (omitting the dependence on $\lambda$ for notational simplicity), for symplectic homology. We similarly use the notation $CF_\bullet(H)$ when $(W,\lambda)$ is fixed.
\end{remark} 


\subsection{Degenerate Hamiltonians and local Floer homology}
In case there is a $1$-periodic orbit of $H$ that is degenerate, we perturb $H$ to a non-degenerate Hamiltonian $\widetilde H$ with the same slope as $H$, choose regular Floer data $(\widetilde J,\widetilde H)$, and define
$$
HF_\bullet(\widehat W, \widehat \lambda, H):= HF_\bullet(\widehat W,\widehat \lambda, \widetilde J,\widetilde H).
$$
\begin{lemma}
This is well-defined, i.e.~it is independent of the choice of perturbation, and of $\widetilde J$.
\end{lemma}

Instead of choosing explicit perturbed Hamiltonians, we package them in local Floer homology, which we now review. Suppose $H$ is a Hamiltonian and assume that $x\in \mathcal P(H)$ is isolated\footnote{In general, we can define local Floer homology for an isolated invariant set.}.
We need the following lemma, which we adapt from \cite{CFHW}:
\begin{lemma}
\label{lemma:isolated_neighborhood}
Suppose that $\gamma$ is an isolated $1$-periodic orbit of $X_H$ with an isolating neighborhood $U$.
Then for every neighborhood $V$ of $\gamma$ with $V \subset U$, there is a $C^2$-small perturbation $\widetilde H$ of $H$ with the following properties:
\begin{itemize}
\item All $1$-periodic orbits of $X_{\widetilde H}$ contained in $U$ are already contained in $V$;
\item For a compatible almost complex structure $\widetilde J$, all Floer trajectories contained in $U$ are already contained in $V$.
\end{itemize}
\end{lemma}

Take a $C^2$-small perturbation $\widetilde H$ as in the lemma so that $1$-periodic orbits in $U$ are non-degenerate (via \cite[Thm.\ ~9.1]{SZ}).
As in \cite{CFHW}, we define the local Floer homology $HF^{loc}_\bullet(\gamma,H)$ of $\gamma$ as the homology of the complex $CF^{loc}_\bullet(U,\widetilde H,\widetilde J)$ generated by $1$-periodic orbits of $\widetilde H$, with differential counting Floer solutions lying in $U$. This is well-defined and independent of the isolating neighborhood $U$, and the perturbed Floer data $(\widetilde J, \widetilde H)$. 

We have the following (see e.g.\ formula (3.1) in \cite{GG15}):
\begin{equation}\label{suppHFloc}
    \mbox{supp}HF_\bullet^{loc}(\gamma,H)\subset [\Delta(\gamma)-n,\Delta(\gamma)+n],
\end{equation}
where $\mbox{supp}HF_\bullet^{loc}(\gamma,H)=\{i: HF_i^{loc}(\gamma,H)\neq 0\}$, and $n=\frac{\dim(W)}{2}$. 

\begin{remark}
\label{rem:free_homotopy}
We observe that the perturbation in Lemma~\ref{lemma:isolated_neighborhood} can chosen such that the $1$-periodic orbits of the perturbed Hamiltonian $\widetilde H$ have the same free homotopy class as $\gamma$.
\end{remark}

\subsection{Spectral sequence}
Suppose now that $H$ is a Hamiltonian that is linear at infinity with $\mathrm{slope}(H)\notin \mathrm{spec}(\alpha)$.
We assume furthermore that the $1$-periodic orbits of $H$ are all isolated.
Hence there are finitely many $1$-periodic orbits with finite action spectrum $\mathcal A_H(\mathcal{P}(H))$. We order the action values in a strictly increasing sequence $\{ a_i \}_{i=1}^k$.
Choose a strictly increasing function $f:\N_0\to \R$ such that $f(i) < a_{i+1}<f(i+1)$.

\begin{proposition}
\label{prop:SS}
There is a spectral sequence converging to the Floer homology $HF_\bullet(W,\lambda,H)$, whose $E^1$-page is given by
\[
E^1_{pq}(H)=\bigoplus_{\substack{\gamma \in \mathcal P(H)\\f(p-1)<\mathcal A_H(\gamma)<f(p)}} HF^{loc}_{p+q}(\gamma,H).
\]
\end{proposition}
We won't give a detailed proof here, but refer to Appendix~B of \cite{KvK} for an almost identical setup. The spectral sequence is the spectral sequence associated with the action filtration given by $f$.
\begin{remark}
\label{rem:degenerate_setup}
This description allows us to define Floer homology for Hamiltonians with isolated possibly \emph{degenerate} orbits.
In addition, we can directly use the free homotopy class of a degenerate periodic orbit, since sufficiently small perturbations cannot change this class as we already observed in Remark~\ref{rem:free_homotopy}. This means that we can decompose Floer homology also in this degenerate setting into free homotopy classes.

A difficulty of this degenerate setup is that a single degenerate orbit can be responsible for several generators in Floer homology. 
Formula~\eqref{suppHFloc} can be used to retain some control.

This point of view is not new, and has been used extensively in for instance \cite{G10} and \cite{H11}. The spectral sequence, although not used in \cite{G10} and \cite{H11}, just gives a convenient packaging, and only serves to make some arguments shorter. The idea of perturbation in Floer theory to get statements of degenerate orbits is even older, and was for example already used in \cite{SZ}.
\end{remark}

\subsection{Index-definiteness and grading.} 
We shall need the following: 
\begin{lemma}
\label{lemma:SH_spread_out}
Suppose that $SH_\bullet(W,\lambda)$ is infinite dimensional, and assume that $\lambda|_{\partial W}$ is an index-definite contact form.
Then $\# \{ i~|~SH_i(W,\lambda)\neq 0 \}=\infty$.
\end{lemma}
\begin{proof}
To prove this, choose a family $\{H_N\}_N$ of admissible Hamiltonians with increasing slopes such that $H_N$ is independent of $N$ on $W$, and so that $CF_\bullet(H_N)$ injects into $CF_\bullet(H_M)$ for $M>N$. By non-degeneracy, each $CF_\bullet(H_N)$ is finitely generated, so the chain complexes get more generators with increasing $N$ (since $\dim SH_\bullet(W,\lambda)=\infty$). 
By the index-definiteness assumption, these new generators have a degree whose absolute value is strictly increasing if $N$ increases sufficiently. This settles the claim.
\end{proof}

\section{Proof of the Generalized Poincar\'e--Birkhoff Theorem}
\label{sec:proof_gpb}
Let $(W,\lambda)$ be a Liouville domain with completion $(\widehat W,\widehat \lambda)$, $r$ the coordinate in the cylindrical end, $B=\partial W$, $\alpha=\lambda\vert_B$, and $\tau$ a Hamiltonian twist map generated by $H=H_t$. The symplectic form on the cylindrical end is $d(r\alpha)$, so by the Hamiltonian twist condition, we get $h_t:B \rightarrow \mathbb{R}^+$ such that $X_{H_t}\vert_B=h_tR_\alpha$. This means that $H_t\vert_{r=1}\equiv C_t>0$, with $\partial_rH_t\vert_{r=1}=h_t$. The family of Hamiltonians $H_t$ is not necessarily linear at infinity, and might hence be unsuitable to compute symplectic homology.
To deal with this we will construct an extension $\widehat H$ to the cylindrical end of $\widehat W$ that is linear at infinity.
By assumption, we have a time-dependent Hamiltonian $H$ defined on $W$.
In a collar neighborhood $\nu(B)$ of the boundary, we will write $H:(1-\epsilon,1] \times B \times S^1 \to \R$, where $r$ is the collar neighborhood parameter.
We extend $H$ to $\widehat H$ on $\widehat W$ using the following procedure.
First of all, define $H_0(b,t):=H(r=1,b,t)$ and $H_1(b,t):=\frac{\partial H}{\partial r}|_{r=1,b,t}$.
We put the remainder in the function $\frac{(r-1)^2}{2!} H_2$, so in short we define
$$
H_2=(H-(H_0+(r-1)H_1)) \frac{2}{(r-1)^2}
$$
on the collar neighborhood $\nu(B)$. By construction $H_2$ is a smooth function on a halfspace.
The functions $H_0$ and $H_1$ are $r$-independent, so admit obvious extensions to $r>1$, but the function $H_2$ is $r$-dependent, so we will appeal to \cite{S}, which is based on reflection, to extend $H_2$ to $r>1$. We call this extension $\overline H_2$. 
Now choose $\delta_1>\delta_0 >0$ and choose a decreasing cutoff function $\rho$ with $\rho|_{[1,1+\delta_0]}=1$ and $\rho(r)=0$ for $r>1+\delta_1$;
\begin{itemize}
\item put  $\widehat H_2(r,b,t) = \overline H_2(r,b,t) \cdot \rho(r)$;
\item put $\widehat H_0(r,b,t) =C \geq \max_t(C_t)$, $\widehat H_1(r,b,t) =A \geq \max_{t,b}(h_t(b))$ for $r\geq 1+\delta_1$;
\item and put $\widehat H_j(r,b,t) = H_j(b,t) \cdot \rho(r) +(1-\rho(r)\, )\widehat H_j(1+\delta_1,b,t)$, for $j=0,1$.
\end{itemize}
The extension is then defined as
\begin{equation}
\label{eq:Ht}
\widehat H:=\widehat H_0(r,b,t) +(r-1) \widehat H_1(r,b,t)+\frac{(r-1)^2}{2!}\widehat H_2(r,b,t)
.
\end{equation}

By the above, we see that $H_0=C_t$ and $H_1=h_t$, so with our choices, we conclude that $\widehat H = A(r-1)+C$ for large $r$.
The extension $\widehat H$ is therefore linear at infinity, and by perturbing $A$ we can assume that $A\notin \mathrm{spec}(\alpha)$.
The same can be arranged for all iterates $\widehat H^{\# k}$ by possibly changing the slope $A$.
The resulting Hamiltonians are then all linear at infinity, but they may have $1$-periodic orbits that are degenerate. 
If all the degenerate $1$-periodic orbits are isolated, then we can still define the Floer homology $HF_\bullet(\widehat W, \widehat \lambda,\widehat H^{\# k})$ using Remark~\ref{rem:degenerate_setup}.
Let us call a Hamiltonian for which all $1$-periodic orbits are isolated, and that is linear at infinity with slope not in the spectrum \emph{weakly admissible}. 

Since we will focus on return maps, it will be convenient to have some shorthand notation.
Define $\widehat\tau^k:=Fl^{X_{\widehat H^{\# k}}}_1$, and with the above remark in mind we put $HF_\bullet(\widehat{\tau}^k):=HF_\bullet(\widehat W, \widehat \lambda,\widehat H^{\# k})$. 
We summarize this discussion in the following lemma:
\begin{lemma}
\label{lemma:extension__Hamiltonian_linear}
The extended Hamiltonians $\widehat H^{\# k}$ are linear at infinity.
Furthermore, if there is an increasing sequence $\{ k_i \}_i\subset \N$ such each Hamiltonian $\widehat H^{\# k_i}$ is weakly admissible, then we have the following isomorphism 
$$
SH_\bullet(W,\lambda) \cong \varinjlim_{k_i} HF_\bullet(\widehat{\tau}^{k_i}),
$$ 
where $\widehat\tau^k=Fl^{X_{\widehat H^{\# k}}}_1$.
$\hfill \square$
\end{lemma}
For later purposes, we need the explicit form of $X_{\widehat H_t}$. This is given by
\begin{equation}\label{X_Hcylend}
\begin{split}
X_{\widehat H_t}&=\left( \partial_r \widehat H_0 +\widehat H_1+ (r-1)\partial_r \widehat H_1
+(r-1)\widehat H_2 +\frac{(r-1)^2}{2} \partial_r \widehat H_2
 \right)R_\alpha\\
&\phantom{=}+\frac{r-1}{r}\left(X_{\widehat H_1}^\xi+\frac{r-1}{2}X_{\widehat H_2}^\xi
-\left(d\widehat H_1(R_\alpha)+\frac{r-1}{2}d\widehat H_2(R_\alpha)
\right)Y \right).
\end{split}
\end{equation}
Here, $Y=r\partial_r$ is the Liouville vector field, and $X_{h}^\xi \in \xi$ is the $\xi$-component of the contact Hamiltonian vector field $X_{h}=hR_\alpha + X_{h}^\xi$ of a Hamiltonian $h: B\rightarrow \mathbb{R}$, defined implicitly by the equation $d\alpha(X^\xi_{h},\cdot)=-dh\vert_{\xi}$.
Due to our choice of interpolation, the second term will be smaller in $C^0$-norm if we choose $\delta_1$ smaller.
We denote the coefficient of $R_\alpha$ by
$$
F= \partial_r \widehat H_0 +\widehat H_1+ (r-1)\partial_r \widehat H_1
+(r-1)\widehat H_2 +\frac{(r-1)^2}{2} \partial_r \widehat H_2.
$$
\begin{lemma}
\label{lemma:scalingF}
If $\delta_1$ is chosen to be sufficiently small, then $F$ is positive.
\end{lemma}
\begin{proof}
To see this, we note that the first three terms are non-negative, and the second term is at least $\min_{t,b}h_t(b)>0$.
To see that the last two terms can be made sufficiently small, note that $\widehat H_2$ has a bound independent of $\delta_1$, and $\partial_r\widehat H_2$ is bounded by $C_2/\delta_1$, where $C_2$ is independent of $\delta_1$. 
Because this term is multiplied by a factor $(r-1)^2$, which is bounded by $\delta_1^2$, the claim follows. 
\end{proof}
As a result we see that $X_{\widehat H}$ is mostly following the positive Reeb direction if we choose $\delta_1$ sufficiently small.
In the proof of Lemma~\ref{indexpositivitylemma} below we will investigate the linearization of $X_{\widehat H}$, which ideally would require closeness to a reparametrized Reeb flow in $C^1$-norm rather than $C^0$-norm. However, $C^1$-closeness does not hold, but we will perform a finer analysis with additional assumptions, which will allow us to fix $\delta_1$.

Lemma~\ref{lemma:extension__Hamiltonian_linear} allows us to compute symplectic homology with the extended Hamiltonian, but it does, by itself, not give any control over periodic orbits in the extension.
To prove our main theorem, we want to show that all generators of $SH_\bullet(W,\lambda)$ represent periodic points of $\tau$ (i.e.\ lie in $W$). To do so, we need to show that the additional periodic points of $\widehat \tau$ do not contribute to the symplectic homology.
Depending on the situation, we will use a filtration by homotopy classes or a filtration by index. More specifically, for $p \in \mbox{Fix}(\widehat \tau^k)$, consider the loop $\gamma_p(t)= Fl^{X_{\widehat H_t}}_t(p)$. Then:

\begin{itemize}
\item If $\dim W=2$, the free homotopy class of $\gamma_p$ in $\tilde \pi_1(W)$ can be used to see that the additional periodic orbits do not contribute homologically;
\item If $\dim W>2$, the CZ-index and the index-definiteness assumption will be used to arrive at the same conclusion. 
\end{itemize}

\subsection{Filtration by homotopy class}
Assume $\dim W=2$. Let $\mbox{Fix}_\partial(\widehat \tau^k):=\mbox{Fix}(\widehat \tau^k)\cap ([1,+\infty)\times B)$. Given $p\in \mbox{Fix}_\partial(\widehat \tau^k)$, let $[\gamma_p]$ be the free homotopy class in $\tilde \pi_1(\partial_p W)\cong \mathbb{Z}$, where $\partial_p W$ is the connected component of $\partial W$ containing $p$. We denote the absolute value of this integer by $|[\gamma_p]|$.

\begin{lemma}
\label{lem:homotopy_class}
Assume the hypothesis of Thm.\ \ref{thmA}, and that $\dim W=2$. Then there is $A>0$, independent of $k$, such that for all $p\in \mbox{Fix}_\partial(\widehat \tau^k)$, we have $|[\gamma_p]|\geq Ak$.
\end{lemma}
\begin{proof} On each circle component of $B$, choose an angular coordinate $\phi$ such that $R_\alpha =\partial_\phi$. From Eq.\ \eqref{X_Hcylend} and Lemma \ref{lemma:scalingF} we see that the component of $X_{\widehat H}$ in the $\partial_\phi$-direction is bounded from below by some constant $A>0$, e.g.\ $A=\inf_B F$. Iterating, we get a bound of the claimed form $Ak$.
\end{proof}

\begin{corollary}
\label{cor:homotopy_filt}
Suppose $W$ and $\tau$ are as in the assumptions of Thm.\ \ref{thmA}, with $\dim W =2$. Then Thm.\ \ref{thmA} holds.
\end{corollary}

\begin{proof}
To prove the statement, we will argue by contradiction, so we assume that the periods of $\tau$ are bounded: denote the minimal periods by $m_0=1< m_1<\ldots< m_M$; we include $m_0$ even if $\tau$ has no fixed points.

Fix a positive integer $N$ and let $A$ be as in Lemma~\ref{lem:homotopy_class}.
Let $\delta$ denote a free homotopy class in $\tilde \pi_1(W)$ that is represented by a simple Reeb orbit (a boundary parallel simple loop). For $i\in \{ 1,\ldots,N\}$ and the iterate $i\delta$, from Cor.\ \ref{cor:SH_surface}, we have $\rk SH^{i\delta}_\bullet(W) = 2$ (here we use the notation from App.\ \ref{app:SH_surfaces}).

We now use the assumption that all fixed points of $\tau$ are isolated, and choose $k>m_M$ such that $k$ is not divisible by $m_1,\ldots,m_M$ (for example choose a large prime).
This choice of $k$ forces the 1-periodic orbits of $\widehat H^{\#k}$ to be isolated on the interior of $W$.
By construction, the Hamiltonian $\widehat H^{\#k}$ is linear at infinity, so we find $r_\infty$ such that $\widehat H^{\#k}$ is linear on $[r_\infty,\infty)\times B$.
If there are non-isolated $1$-periodic points of $\widehat H^{\#k}$ on the cylindrical part $[1,r_\infty]\times B$, then we use Lemma~\ref{lemma:perturb_cylinder} below to perturb the Hamiltonian $\widehat H^{\#k}$; this perturbation makes all orbits on $[1,r_\infty]\times B$ non-degenerate, and does not affect $1$-periodic orbits on the interior of $W$.
Hence we obtain a weakly admissible, possibly degenerate Hamiltonian, which we continue to denote by $\widehat H^{\#k}$.
For this Hamiltonian we can define Floer homology using Remark~\ref{rem:degenerate_setup}.

By choosing an increasing sequence $\{ k \}$ of primes, we can then define $SH_\bullet^{i\delta}(W) =\varinjlim_k HF_\bullet^{i\delta}(\widehat H^{\# k})$.
Hence we find a sufficiently large $k$ that
$$
2\leq \rk HF_\bullet^{i\delta}(\widehat H^{\# k}) 
= \sum_{p,q}E^\infty_{pq}(\widehat H^{\# k})
\leq \sum_{p,q}E^1_{pq}(\widehat H^{\# k})
=\sum_{p,q}\rk HF_{p+q}^{loc,i\delta}(\gamma,\widehat H^{\# k})
.
$$
All of these sums are finite by the assumption that the fixed points are isolated.
We conclude that there is a $1$-periodic orbit $\gamma_{k,i\delta}$ of $\widehat H^{\# k}$ whose free homotopy class equals $i\delta$.
From Lemma~\ref{lem:homotopy_class}, every $p\in \mbox{Fix}_\partial(\widehat \tau^k)$ has $ [\gamma_p]=j \delta$ with $j\geq Ak$. If we choose $k>N/A$ we see that $j >N$, so the $1$-periodic orbit $\gamma_{k,i\delta}$ is represented by a fixed point of $\tau^k$.

This works for all $N$, so by sending $k$ to infinity we get infinitely many periodic points of $\tau$. To see that these are geometrically distinct, note that if $p\in \mbox{Fix}(\tau^k)$ and $a:=[\gamma_p]=i\delta$ is boundary parallel, then $\gamma_p^\ell$ has homotopy class $\ell a$, so another orbit must represent the free homotopy class $a$. Taking the limit in $k$, we see that new generators in the homotopy class $a$ need to appear in order to generate $SH^a(W)$. This gives infinitely many geometrically distinct interior periodic points (in different boundary parallel homotopy classes).
\end{proof}

\subsection{Filtration by index}
We now deal with the second case, so we assume now that $\dim W>2$, $c_1(W)|_{\pi_2(W)}=0$, and that the Reeb flow is strongly index-definite. To set up the argument, we first need to establish that index-definiteness of the linearized Reeb flow equation at the boundary (in the sense of Definition~\ref{def:indexposODE} in Appendix~\ref{app:indexposODEs}) implies index-definiteness of the linearized Hamiltonian equation along the cylindrical end:

\begin{lemma}\label{indexpositivitylemma}
Assume that $(\xi|_{B},d\alpha|_B)$ is symplectically trivial, and that the linearized Reeb flow equation $\dot \psi = \nabla_\psi R_\alpha$ along $B=\partial W$ is strongly index-definite. Then, the linearized Hamiltonian flow equation $\dot \psi=\nabla_\psi X_{\widehat H} $ of the extension $\widehat H$ given by Equation~\eqref{eq:Ht} is also strongly index-definite along the cylindrical end $[1,+\infty)\times B$.
\end{lemma}

\begin{proof} 
We prove this using a matrix representation.
To do this, we need to symplectically trivialize the full tangent bundle on the cylindrical ends.
Given a symplectic trivialization of $(\xi|_{B},d\alpha|_B)$, we only need to trivialize the symplectic complement of $\xi$.
We do this using the trivialization $L=\langle Y=r \partial_r, R \rangle$, where $R=R_\alpha/r$ is the Reeb vector field at the $r$-slice.

We will work with the usual formalism of time-dependent Hamiltonians, and we do not include this time-dependence in the notation. Exterior and covariant derivatives are computed using the base manifold only, and do not involve time derivatives. 
We will also use the following notation:
$$
X_\xi:=X_{\widehat H_1}^\xi+\frac{r-1}{2}X_{\widehat H_2}^\xi ,
$$
$$
G:=d\widehat H_1(R_\alpha)+\frac{r-1}{2}d\widehat H_2 (R_\alpha).
$$
To compute the linearization, we choose a convenient connection $\nabla$, namely the Levi-Civita connection for the metric $1/r^2 \cdot dr\otimes dr +\alpha\otimes \alpha +d\alpha(\cdot,J\cdot)$.
This connection has the following properties:
\begin{itemize}
\item $\nabla Y=0$. Keep in mind that $Y$ is the Liouville vector field $r\partial_r$;
\item $\nabla_{R_\alpha} R_\alpha=0$ and $\nabla_Y R_\alpha=0$;
\item $\nabla_X R_\alpha \in \xi$ for all $X\in \xi$.
\end{itemize}
With respect to this connection we compute the linearization as
\begin{equation}
\label{linearizedXH}
\nabla X_{\widehat H}= F \nabla R_\alpha + dF \otimes R_\alpha + \frac{1}{r^2}dr \otimes (X_\xi-GY)+\frac{r-1}{r}(\nabla X_\xi-dG \otimes Y).
\end{equation}
Before we continue our analysis of the linearization, we first need to discuss the behaviour of the Hamiltonians $\widehat H_j$ and their derivatives under rescaling the interpolation parameter $\delta_1$.
We will write the terms in the expression~\eqref{eq:Ht} as $\widehat H_j'$ if we use $\delta_1'$ as interpolation parameter.
For $\delta_1'<\delta_1$ we have the following:
\begin{itemize}
\item derivatives in the $B$-direction (denoted $\partial_b$) admit a uniform bound, independent of $\delta_1$, i.e.
$$
\max_{[1,+\infty)\times B} |\partial^k_b \widehat H_j'| \leq \max_{[1,+\infty)\times B} |\partial^k_b \widehat H_j| \mbox{ for all } k\geq 0; 
$$
\item derivatives in the $r$-direction scale as follows:
$$
\max_{[1,+\infty)\times B} |\partial_r^k \widehat H_j'| \leq 
\left(\frac{\delta_1}{\delta_1'}\right)^k \max_{[1,+\infty)\times B} |\partial_r^k \widehat H_j| \mbox{ for all } k \geq 0. 
$$
\end{itemize}
Keeping this scaling behaviour in mind, we regroup terms in Equation~\eqref{linearizedXH} to obtain the following representation:
\[
\nabla X_{\widehat H}
=L_0 + L_1,
\]
where
$$
L_0=F \nabla R_\alpha + dF \otimes R_\alpha +\frac{1}{r^2}dr \otimes (X_\xi-GY)
-\frac{r-1}{r^2} dG(Y) dr \otimes Y
+\frac{r-1}{r^2} dr \otimes \nabla_Y X_\xi
$$
and
$$
L_1=\frac{r-1}{r}\left(
\nabla^\xi X_\xi +\alpha \otimes \nabla_{R_\alpha} X_\xi
-R_\alpha(G) \alpha \otimes Y-d^\xi G \otimes Y
\right).
$$
Here, $\nabla^\xi= P_\xi \nabla\vert_\xi,$ where $P_\xi$ is the orthogonal projection to $\xi$, and $d^\xi=d\vert_\xi$. We will explain below that the matrices $L_0$ and $L_1$ have the following matrix representations, with respect to the decomposition $T\widehat W=\xi \oplus \langle Y,R\rangle$:
\[
L_0=
\left(\begin{array}{@{}c|c@{}}
  F \cdot \nabla^\xi R_\alpha
  & \begin{matrix}
  U & 0 \\
  V & 0 \\
  \end{matrix} \\
\hline
  \begin{matrix}
  0 & 0\\
  W & Z
  \end{matrix} &
  \begin{matrix}
     a & 0 \\
     b & c
  \end{matrix}
\end{array}\right), 
\quad
L_1=\frac{r-1}{r}
\left(\begin{array}{@{}c|c@{}}
  \nabla^\xi X_\xi
  & \begin{matrix}
  0 & U' \\
  0 & V' \\
  \end{matrix} \\
\hline
  \begin{matrix}
  W' & Z'\\
  0 & 0
  \end{matrix} &
  \begin{matrix}
     0 & a' \\
     0 & 0
  \end{matrix}
\end{array}\right)
\]
This is clear for $L_0$. We further want to show that $L_0\in \mathfrak{sp}(2n)$, which will constrain the entries more. 
Since we know that $L_0+L_1\in \mathfrak{sp}(2n)$, we will show that $L_1\in \mathfrak{sp}(2n)$, since the latter contains fewer terms.
For this we note the following:
\begin{itemize}
\item the matrix representation for $\nabla^\xi X_\xi$ is in $\mathfrak{sp}(2n-2)$. This is because these entries come from the $\xi$-part of a contact Hamiltonian;
\item the matrix representation for $R_\alpha(G) \alpha \otimes Y$ is in $\mathfrak{sp}(2)$: the non-trivial entry corresponds to the element $a'$;
\item non-trivial entries in the matrix representation of $-d^\xi G \otimes Y$ appear only on the first row of the lower left block.
These correspond to the elements $W',Z'$;
\item non-trivial entries in the matrix representation of $\alpha \otimes \nabla_{R_\alpha} X_\xi$ appear as the last column. 
We will show that these correspond to the elements $U'$ and $V'$.
We claim that $\langle \nabla_{R_\alpha} X_\xi ,R_\alpha \rangle=0$. Indeed, since the contact structure is orthogonal to the Reeb vector field with our choice of metric, we have
$$
0= R_\alpha
\langle X_\xi ,R_\alpha \rangle
=\langle \nabla_{R_\alpha} X_\xi ,R_\alpha \rangle
+\langle X_\xi ,\nabla_{R_\alpha} R_\alpha \rangle
=\langle \nabla_{R_\alpha} X_\xi ,R_\alpha \rangle
.
$$
Similarly, we obtain $\langle \nabla_{R_\alpha}X_\xi,Y\rangle=0$. This means that the $L$-entries in the matrix representation of $\alpha \otimes \nabla_{R_\alpha} X_\xi$ are zero.
\item we have $(W^\prime,Z^\prime)^T=J\cdot (U^\prime, V^\prime)^T=(-V^\prime,U^\prime)^T.$ This follows since $-d^\xi G$ is dual to $\nabla_{R_\alpha}X_\xi$, i.e.\ $d\alpha(\nabla_{R_\alpha}X_\xi,\cdot)=-d^\xi G$.
\end{itemize}
We conclude that $L_1\in \mathfrak{sp}(2n)$, and hence $L_0$ is, too.
Observe also that for all $\epsilon>0$ we can choose $\delta_1>0$ such that $\Vert L_1 \Vert <\epsilon$ due to the scaling behaviour we discussed earlier: this can be done in a way that is compatible with Lemma~\ref{lemma:scalingF}, i.e.~$\delta_1$ getting smaller as $\epsilon$ gets smaller.

Since $J_0L_0$ is symmetric, we can fix the terms of $L_0$. They must necessarily have the following form:
\[
L_0=
\left(\begin{array}{@{}c|c@{}}
  F \cdot \nabla^\xi R_\alpha
  & \begin{matrix}
  U & 0 \\
  V & 0 \\
  \end{matrix} \\
\hline
  \begin{matrix}
  0 & 0\\
  V & -U
  \end{matrix} &
  \begin{matrix}
     a & 0 \\
     b & -a
  \end{matrix}
\end{array}\right) \in \mathfrak{sp}(2n).
\]
This matrix has precisely the form that we consider in Appendix \ref{app:indexposODEs}. Moreover, note that strong index-definiteness is invariant under scaling by a positive (possibly time-dependent) function of the generating matrix. Indeed, this scaling has the effect of positively reparametrizing the flow, and so the new flow intersects the Maslov cycle as often as the original one (although the constants in the definition of strong index-definiteness might change). 
Therefore, since the ODE $\dot \psi=\nabla^\xi_\psi R_\alpha$ is strongly index-definite by assumption and $F>0$, then so is the ODE $\dot \psi=F\cdot (\nabla^\xi_\psi R_\alpha)$. 
Lemma \ref{App:indexposestimate} in Appendix~\ref{app:indexposODEs} now tells us that the system
$\dot \psi =L_0 \psi$ is strongly index-definite.
By choosing $\delta_1$ sufficiently small, we can make the matrix $L_0$ get arbitrarily $C^0$-close to $L_0+L_1= \nabla X_{\widehat H}$.
Since the system $\dot \psi =L_0 \psi$ is strongly index-definite, we can adapt Lemma~2.2.9 from \cite{U99} to see that $\dot \psi = \nabla_\psi X_{\widehat H}$ is strongly index-definite, too.
This concludes the proof of Lemma \ref{indexpositivitylemma}. 
\end{proof}

We need the following lemma to ensure that our Hamiltonians are weakly admissible.
\begin{lemma}
\label{lemma:perturb_cylinder}
Given an extension $\widehat \tau: \widehat{W} \to \widehat{W}$ as in the beginning of Section~\ref{sec:proof_gpb}, there is a Hamiltonian perturbation $\widetilde \tau=\phi_f^1 \circ \widehat \tau: \widehat{W} \to \widehat{W}$ with the following properties:
\begin{enumerate}
\item $\phi_f^1(x) =x$ for all $x$ not in a neighborhood of $[1,r_\infty]\times B$, for some fixed $r_\infty>1$. In particular, all interior fixed points of $\tau$ are unaffected by the perturbation;
\item all fixed points of $\widetilde \tau|_{[1,r_\infty]\times B}$ are non-degenerate and hence isolated;
\item by the standard composition rule for Hamiltonians $\widetilde \tau$ is the time $1$-flow of a Hamiltonian $\widetilde H$. This Hamiltonian $\widetilde H$ is $C^2$-close to $\widehat H$, and its fixed points have Robbin-Salamon index close to that of the unperturbed fixed points.
\end{enumerate}
\end{lemma}

\begin{proof}
We adapt the argument from \cite[Lemma~2]{B} to our setting.
Set $V:=[1,r_\infty]\times B$, where $r_\infty>1$ is such that $\widehat H$ is linear on $[r_\infty,\infty) \times B$. 
We need to find a $C^2$-small function $f$ vanishing on the complement of a neighborhood of $V$ such that all fixed points of $\phi_f^1 \circ \widehat \tau$ are non-degenerate on $V$.
Consider the map
$$
j:V \longrightarrow V \times V,\quad
x\longmapsto (x, \widehat \tau(x)), 
$$
and denote its image by $\Gamma$. Define the diagonal $\Delta = \{ (v,v) \in V\times V~|~v\in V \}$.
Observe that $v\in \mbox{Fix}(\widehat \tau)$ is a non-degenerate fixed point of $\widehat \tau$ if and only if $\Gamma$ and $\Delta$ intersect transversely at $(v,v)$.

For all points $(v,v)\in \Gamma \cap \Delta$ choose a Darboux ball $B_{\epsilon}(v) \subset \widehat{W}$ such that $v$ corresponds to $0$ in the Darboux ball.
To choose $\epsilon$, let $\lambda_{max}$ denote the maximal (in absolute value) eigenvalue of $d_x\widehat \tau$ over all fixed points of $\widehat \tau$ in the cylinder. In a formula,
$$
\lambda_{max}:=
\max\{
\vert \lambda \vert ~|~\lambda \text{ eigenvalue of } d_x\widehat \tau, x\in \mbox{Fix}(\widehat \tau) \cap [1,r_\infty] \times B \}.
$$ 
Note that $\lambda_{max}\geq 1$, since $d_x\widehat \tau$ is symplectic. Choose $\epsilon$ so small such that the following two properties hold:
\begin{itemize}
\item if $x\in B_{\epsilon/2\lambda_{max}}(v)$, then $\widehat \tau(x) \in B_{3\epsilon/4}(v)$. 
\item for all interior fixed points of $\tau$, i.e.~for $x\in \mbox{Fix}(\tau|_{int(W)})$ we have $d(x,[1,\infty]\times B)>\epsilon$, where $d$ is some fixed reference metric (for example, induced by the Riemannian metric $\widehat \omega(\cdot,\widehat J\cdot)$).
\end{itemize}
We give some intuition for these choices, before going into the computation.
By the first property, we retain some control after applying $\widehat \tau$ to a point that is sufficiently close to a fixed point.
Intuitively, if $x$ is close to the fixed point $v$, then by the definition of $\lambda_{max}$, the map $\widehat \tau$ sends $x$ approximately away by a factor of at most $\lambda_{max}$, and so we ensure that if $x\in B_{\epsilon/2\lambda_{max}}(v)$, then $\widehat \tau(x)\in B_{3\epsilon/4}(v)$.
Below we will define Hamiltonian functions to perturb the map $\widehat \tau$, and this property will ensure that we have maximal control over the value of the Hamiltonian vector fields.
This point is actually not essential, but it makes the computation below a little more uniform.

We now come to our Hamiltonian perturbation functions.
Choose functions $f_{v,i}$ for $i=1,\ldots,2n$ such that, in Darboux coordinates $z=(z_1,\dots,z_{2n})$, we have $f_{v,i}(z) = z_i\cdot \rho_v(z)$, and $\rho_v$ is a cutoff function that equals $1$ on $B_{3\epsilon/4}(v)$ and vanishes outside $B_{\epsilon}(v)$. For the sake of explicitness, note that the Hamiltonian vector field of $f_{v,i}$ is given by
$$
X_{f_{v,i}} = \rho_v(z) J_0 \cdot e_i+ z_i \cdot X_{\rho_v}, 
$$
where $e_i$ is the $i$-th standard basis vector, and $J_0$ is the standard complex structure on the Darboux ball $B_{\epsilon}(v)$. 
By construction, this vector field vanishes on the complement of $B_{\epsilon}(v)$. Moreover, for sufficiently small $r$, the time-$1$ flow of the Hamiltonian vector field of $rf_{v,i}$ on the smaller ball $B_{\epsilon/2\lambda_{max}}(v)$ is the map $z\mapsto z+r J_0 \cdot e_i$. If $x\in \mbox{Fix}(\widehat \tau)\cap B_{\epsilon/2\lambda_{max}}(v)$ we have
$$
\frac{\partial}{\partial r}\Big\vert_{r=0}
\phi^1_{r f_{v,i}} \circ \widehat \tau(x)
=\frac{\partial}{\partial r} \Big\vert_{r=0}
\phi^r_{f_{v,i}} \circ \widehat \tau(x)
=X_{f_{v,i}}(\widehat \tau(x)).
$$ 
For completeness, we observe that each of these functions $f_{v,i}$ is a $C^\infty$ function defined on all of $\widehat W$, vanishing outside $B_\epsilon(v)$.

Since $\Gamma \cap \Delta$ is compact, we find a finite cover of its projection to $V$ of the form $\bigcup_{v \in I} B_{\epsilon/2\lambda_{max}}(v)$.
We make the following observation. Consider $x\in \mbox{Fix}(\widehat \tau)\cap [1,r_\infty]\times B$. 
Then there is $v\in I$ such that $x\in B_{\epsilon/2\lambda_{max}}(v)$.
On this small ball, the Hamiltonian vector fields associated with $f_{v,1},\ldots, f_{v,2n}$ are linearly independent, and form a basis of sections.

Define the finite-dimensional vector space
$$
R:= \left(\R^{2n} \right)^{\# I } =\R^{D},  
$$
where we have set $D:=2n\# I$. We relabel the functions $f_{v,i}$ using a single index $j$, and put $f=(f_1,\cdots,f_D)$. Define the projection $p:V \times V \times R \to R,~(x,y,r)\mapsto r$, and consider the ``universal'' space 
$$
\Gamma_R=\left\{ (x,y,r) \in V \times V \times R ~\Bigg|~x\in V,r=(r_1\ldots,r_D)\in R,\;y=\phi_{r\cdot f}^1 \circ \widehat \tau(x), \; r\cdot f=\sum_j r_j f_{j}\right\}.
$$ 
Note that the function $r\cdot f$ is a $C^\infty$ function defined on all of $\widehat W$: this function vanishes outside a neighborhood of $[1,r_\infty] \times B$.

\medskip

\noindent
{\bf Claim: } The space $\Gamma_R$ intersects the enlarged diagonal $\Delta_R = \{ (v,v,r)\in V\times V \times R \} $ transversely for $r$ that are sufficiently close to $0$.
In particular, $V_R:=\Gamma_R \cap \Delta_R$ is a submanifold.

\begin{proof}[Proof of claim]
To verify the claim, we compute the derivatives of the map 
$$
j_R:(x,r)\mapsto (x,\phi_{r\cdot f}^1 \circ \widehat \tau(x),r)
$$ 
and the corresponding map for $\Delta_R$, $j_\Delta: (x,r)\mapsto (x,x,r)$. 
For $j_R$ we find the derivative
$$
d_{x,r=0} j_R=
\left(
\begin{array}{cc}
\id_V & 0 \\
d_x \widehat \tau & X_{f_1}(\widehat\tau(x)),\ldots, X_{f_D}(\widehat\tau(x)) \\
0 & \id_R
\end{array}
\right).
$$
For $j_\Delta$ we find the derivative
$$
d_{x,r=0} j_\Delta
=
\left(
\begin{array}{cc}
\id_V  & 0\\
\id_V & 0 \\
0 & \id_R
\end{array}
\right).
$$
Given a point $(x,x,0) \in V_R=\Gamma_R \cap \Delta_R$ (so $\widehat \tau(x)=x$), there is $v\in I$ such that $x\in B_{\epsilon/2\lambda_{max}} (v)$.
By construction, the vector fields $X_{f_{v,1}}\circ \widehat \tau,\ldots, X_{f_{v,2n}}\circ \widehat \tau$ are linearly independent on $B_{\epsilon/2\lambda_{max}} (v)$.
This means that, taken together, the matrix representations of $d_{x,r=0} j_R$ and $d_{x,r=0} j_\Delta$ have full rank, namely $2n+2n+D$.
We conclude that $j_R$ is transverse to the enlarged diagonal $\Delta_R$ for $r=0$, and hence, by compactness, also for small $r$.
\end{proof}
Applying Sard's theorem to the projection $p|_{V_R}$, we find a regular value $r_0$ of $p|_{V_R}$ close to $0$.
We see that
$$
\Gamma_{r_0}=\{ (x,y) \in V \times V~|~y=\phi_{r_{0}\cdot f}^1 \circ \widehat \tau(x) \} 
$$
intersects $\Delta$ transversely.
This means that all fixed points of $\phi_{f_0}^1 \circ \widehat \tau$ in $V$ are non-degenerate, where $f_0=r_0\cdot f$, so claim (2) holds. Since we can choose the regular value $r_0$ arbitrarily small, and since the support of the perturbation $f_0$ is a small neighborhood of $[1,r_\infty] \times B$, the claim (1) holds. To see that claim (3) also holds, we note that the free homotopy class of a $1$-periodic orbit is not affected by this perturbation if it is sufficiently small. For the index, we use the same argument as before. The unperturbed system is strongly index-definite, and the same will be true for small perturbations. This concludes the proof of the lemma.
\end{proof}

\begin{proof}[Proof of Thm.~A ($\dim W>2$)]
Write $\tau =\phi^1_H$ for $H$ as in Def.\ \ref{def:twistmap}. Assuming its interior fixed points are isolated, we have finitely many isolated interior $1$-periodic orbits of $H$, say $\gamma_1,\ldots, \gamma_k$. The starting points $\gamma_1(0),\ldots,\gamma_k(0)$ are the fixed points of $\tau$.

Assume by contradiction that the minimal periods of all interior periodic points of $\tau$ are, in increasing order, given by $m_0=1,m_1,\ldots, m_\ell$.
Take an increasing sequence $\{ p_i\}_{i=1}^\infty$ going to infinity, and such that each $p_i$ is indivisible by the $m_1,\ldots, m_\ell$. 
For instance, one can take the sequence $\{ p_i\}$ to be an increasing sequence of primes all of which are larger than $\max_{j} m_j$.

As in the proof of Corollary~\ref{cor:homotopy_filt}, we can appeal to Remark~\ref{rem:degenerate_setup} to define Floer homology for a possibly degenerate Hamiltonian.
Indeed, due to the choice of $p_i$'s, all fixed points of $\widehat \tau^{p_i}$ are isolated, and we can apply Lemma~\ref{lemma:perturb_cylinder} if necessary to perturb the Hamiltonian $\widehat H^{\# p_i}$ on the cylindrical part $[1,r_\infty]\times B$, for some $r_\infty$. This ensures that the Hamiltonian $\widehat H^{\# p_i}$ is weakly admissible, so we can use local Floer homology and the spectral sequence from Proposition~\ref{prop:SS}, to define $HF(\widehat H^{\# p_i})$.
Hence we can compute symplectic homology as 
$
SH_\bullet(W)=\varinjlim_i HF_\bullet(\widehat H^{\#p_i})$. 
By Lemma~\ref{lemma:SH_spread_out}, for all $N>2n k$, where $\dim(W)=2n$, we find distinct degrees $i_1,\ldots,i_N$ such that $SH_{i_j}(W)\neq 0$, ordered by increasing absolute value. 
By Lemma~\ref{indexpositivitylemma}, we can choose $p_i$ sufficiently large such that the following hold:
\begin{enumerate}
\item Each $1$-periodic orbit of $\widehat{H}^{\# p_i}$ that is contained in $\widehat W \setminus \mbox{int}(W)$ has RS-index whose absolute value is larger than $|i_N|+2n$;
\item the Floer homology groups $HF_{i_j}(\widehat H^{\# p_i})$ are non-trivial for $j=1,\ldots, N$. 
\end{enumerate}
Now consider the spectral sequence from Proposition~\ref{prop:SS} for $\widehat H^{\# p_i}$.
We deduce from (2) that there must be non-trivial summands on $E^1_{pq}(\widehat H^{\# p_i})$ with $p+q=i_j$ for $j=1,\ldots, N$.
Since the terms of the spectral sequence are made up from local Floer homology groups, and we know from (1) that no $1$-periodic orbit in $\widehat W \setminus \mbox{int}(W)$ can contribute to local Floer homology of degree $i_j$, we conclude that every term $E^1_{pq}(\widehat H^{\# p_i})$ in the spectral sequence with $p+q=i_j$ must come from the local Floer homology of an orbit $\gamma$ in int$(W)$.

Because we have assumed that the $p_i$'s are indivisible by $m_1,\ldots, m_\ell$ we conclude that each such orbit $\gamma$ must be an iterate of one of the orbits $\gamma_1,\ldots,\gamma_k$. Moreover,  by~\eqref{suppHFloc} and Section \ref{sec:def_indices}: $$\mbox{supp}HF^{loc}_\bullet(\gamma_j^{p_i},\widehat H^{\# p_i})\subset[p_i \Delta(\gamma_j)-n, p_i \Delta(\gamma_j)+n].$$ This covers at most $2n k$ different degrees, leaving some of the degree $i_j$ uncovered as we had chosen $N>2n k$. This is a contradiction.
\end{proof}

\begin{proof}[Proof of Thm~\ref{thmB}]
We only need to show that $\dim SH_\bullet(W)=\infty$. Since $W\subset T^*M$ is star-shaped, from Viterbo's theorem \cite{V99} we have $SH_\bullet(W)\cong H_{\bullet}(\mathcal{L} M;\Z_2)$ where $\mathcal{L}M$ is the free loop space of $M$. The statement is more subtle when using $\Z$ or $\mathbb{Q}$-coefficients, see \cite{A15}. Now we can apply the following theorem due to Gromov: 
\begin{thm}\cite[Sec.\ 1.4]{G78}\label{Gromov}
Let $(M,g)$ be a closed Riemannian manifold with finite fundamental group. For $a>0$, let $\mathcal{L}M$ be the free loop space of $M$, and let $\mathcal{L}^{<a}M\subset \mathcal{L}M$ denote the space of free loops with length less than $a$. Let $\iota^a: \mathcal{L}^{<a}M \hookrightarrow \mathcal{L}M$ denote the inclusion, and $\iota^a_k: H_k(\mathcal{L}^{<a}M;\mathbb{R})\rightarrow H_k(\mathcal{L}M;\mathbb{R})$ be the map induced in real homology of degree $k$. Then there exists a positive constant $C=C(M,g)$ such that
$$
\sum_{k\geq 0}\mbox{\text{rank}} (\iota^a_k)\geq Ca.
$$
\end{thm}
Together with the above, this tells us that $SH_\bullet(W)$ is infinite dimensional.
\end{proof}

\appendix

\section{Hamiltonian twist maps: examples and non-examples}
\label{sec:examples} We will now discuss some examples that help clarify the nature of the Hamiltonian twist condition. 

\subsection{Examples.} The following construction, an adaptation of a standard one, further illustrates that the Hamiltonian twist condition is \emph{not} localized at $B$.

\begin{proposition}
For each $k,\ell\in \N$, there are strict contact manifolds $(Y_k,\alpha_{k,\ell})$ carrying adapted open books $(B_k= B,\pi_k)$, $\pi_k: Y_k\backslash B \rightarrow S^1$, with fixed page $\Sigma$, such that the following holds:
\begin{itemize}
\item The return maps $\tau_k$ all agree in a collar neighborhood of $B=\partial \Sigma$, and are generated by Hamiltonians $H_k$;
\item Furthermore, there is a symplectomorphism $\phi_k$ from $\mathring \Sigma$, the interior of the page $\Sigma$, to the open subset $W_2$ of the Liouville completion $\widehat W$ of a fixed Liouville domain $(W,\lambda)$ with $\partial W=B$, where
$$
W_2= W \cup_\partial ([1,2)\times B, d (r\alpha_B)\, ),
$$
and $\alpha_B=\lambda\vert_B$ is the contact form at $B$, and $r \in [1,2)$.
\item The return map $\phi_k\circ \tau_k \circ \phi_k^{-1}$ extends to a Hamiltonian diffeomorphism $\bar \tau_k$ on the closure $\bar W_2$, generated by Hamiltonians $\bar H_k$.
\item The Hamiltonian twist condition holds for $\bar H_k$ for $k\leq \ell$, but not for $k>\ell$.
\end{itemize}
\end{proposition}
\begin{proof}
Consider a Liouville domain $(W,\lambda)$ with a $2\pi$-periodic Reeb flow on its boundary (e.g.\ $D^*S^2$).
We identify a collar neighborhood $\nu_W(B)$ of $B=\partial W$ with $(1/2,1]\times B$, where $B=\{r=1\}$, via a diffeomorphism $\varepsilon: (1/2,1]\times B \longrightarrow \nu_W(B)\subset W$. We assume $\lambda=r\alpha_B$ along $\nu_W(B)$, $\alpha_B=\lambda\vert_B$. Define the smooth Hamiltonian
\[
H(x)=\begin{cases}
0, & \text{if }x\notin \nu_W(B), \\
f(r), & \text{if }x=\varepsilon(r,b) \in \nu_W(B).
\end{cases}
\]
Here $f$ is a smooth, decreasing function with the property 
\begin{itemize}
\item $f(1/2)=0$;
\item $f'(r)\geq -2\pi$ and $f'(r)=-2\pi$ near $r=1$.
\end{itemize}
The Hamiltonian vector field of $H$ is given by
\[
X_H(x)=\begin{cases}
0 & \text{if }x\notin \nu_W(B), \\
f'(r) R_\alpha & \text{if }x=\varepsilon(r,b) \in \nu_W(B).
\end{cases}
\]
Define the fibered Dehn twist by $\tau(x)=Fl^{X_H}_1(x)$, where $Fl^{X_H}_t$ is the Hamiltonian flow of $H$ with respect to $d\lambda$. We have $\tau^*\lambda=\lambda - dU$, where we choose the primitive $U$ to be a negative function: with a computation we can show that it is possible to choose $U(1)=-2\pi$, and will do so. The iterate $\tau^k$ is generated by $H_k=kH$, and $(\tau^k)^*\lambda=\lambda-dU_k$, with $U_k=\sum_{j=0}^{k-1}(U \circ \tau^j)$.

We consider the associated open book
$$
Y_k=OB(W,\tau^k):= B\times D^2 \cup_\partial W_{\tau^k},
$$
where $W_{\tau^k}=W\times \R /(x,t)\sim(\tau^k(x),t+U_k(x))$ is the mapping torus. 
The manifold $Y_k$ carries an adapted contact form $\alpha_{k,\ell}$ which looks like $\alpha_{k,\ell}=\lambda+d\theta$ along $W_{\tau^k}$, and $\alpha_{k,\ell}=h_1(\rho)\alpha_B+h_2(\rho)d\theta$ along $B\times D^2$. Here, $(\rho,\theta) \in D^2$, and $h_1$ and $h_2=h_{2,k,\ell}$ are suitable profile functions, which we will fix now. Choose $h_1$ and $h_2$ such that:
\begin{itemize}
\item they do not depend on $k$ for $\rho\leq 1/2$;
\item $h_1^\prime\leq 0$ with equality only at $\rho=0$. We may take $h_1(\rho)=2-\rho^2$ near $\rho=0$ (this is not essential but very convenient);
\item near $\rho=0$ we have $-\frac{h_2^\prime}{h_1^\prime}(\rho)=\ell + \epsilon>0$ (non-singular) for some small $\epsilon \in (0,1)$.
\item $h_2\equiv k$, $h_1=-\rho+2$ near $\rho=1$ (so $h_2$ depends on $k$ on the interval $(1/2,1]$).
\end{itemize}
Note that, in the definition of $Y_k$, the binding model is glued to the mapping torus using the gluing map
\[
\begin{split}
\Phi_{glue}: B\times D^2_{\rho>1/2} & \longrightarrow W_{\tau^k}\\
(b;\rho,\theta) & \longmapsto \left(2-\rho,b;\frac{-U_k(1)\theta}{2\pi}\right)= (2-\rho,b;k \theta).
\end{split}
\]
This pulls back $d\theta+ \lambda$ to $kd\theta +(2-\rho) \alpha_B$. This explains the above choices.

The global hypersurface of section, i.e.\ a fixed page, is $\Sigma=W \cup_\partial B \times [0,1]$, with coordinate $\rho\in [0,1]$, and we can compute the return map $\tau_k$ explicitly. We find:
\[
\tau_k(x)=\begin{cases}
\tau^k(x), & \text{if }x \in W, \\
(Fl^R_{-2\pi h_2'(\rho)/h_1'(\rho)}(b),\rho), & \text{if }x=(b,\rho) \in B\times [0,1],
\end{cases}
\]
where $Fl^R_t$ is the Reeb flow of $\alpha_B$ at $B$. The Hamiltonian generating $\tau_k$ can be obtained by patching $H_k$ on $W$ to a Hamiltonian that generates $\tau_k$ along $B\times [0,1]$; we need to match the slopes on the boundary, which can be done by rewriting $\tau_k(b,\rho)=(Fl^R_{-2\pi (h_2'(\rho)/h_1'(\rho)+k)}(b),\rho)
$. 
Then $H_k$ extends to $\Sigma$ via $\bar H_k(r)=-2\pi\int_1^\rho(h_2'(s)/h_1'(s)+k)h_1^\prime(s)ds + f(1)$ along $B\times [0,1]$. Note that the form $d\lambda$ also extends along $B\times [0,1]$ via $d\alpha_{k,\ell}\vert_\Sigma=h_1^\prime(\rho)d\rho \wedge \alpha_B$. 
Therefore $H_k$ generates $\tau_k$, and $\tau_k$ is independent of $k$ on the collar neighborhood $B\times [0,1/2]$.

To complete the proof, we first note that the $2$-form $d (h_1(\rho) \alpha_B)$ is degenerate on $\partial \Sigma$.
However, the map 
\[
\begin{split}
\phi_k: {\mathring \Sigma} &\longrightarrow W_2= W \cup_\partial ([1,2) \times B, d (r\alpha_B)\, ),\\
w &\longmapsto \begin{cases}
w & w \in W \\
(h_1(\rho),b) & w =(b,\rho)\in B \times (0,1]
\end{cases}
\end{split}
\]
is a symplectomorphism, and the closure of $W_2$ is an actual Liouville domain. Furthermore, due to our explicit choice $h_1(\rho)=2-\rho^2$ near $\rho=0$, we find $\rho=\sqrt{2-r}$, so we can compute the conjugated return map $\phi_k \circ \tau_k \circ \phi_k^{-1}$ near $r=0$ as
$$
\phi_k \circ \tau_k \circ \phi_k^{-1}(r,b)=
(r, Fl^R_{2\pi (\ell+\epsilon)}(b) ).
$$
This map extends to a symplectomorphism $\bar \tau_k:\bar W_2 \to \bar W_2$. 
Here, note that $\phi_k^{-1}$ is not smooth at $r=2$, but this is resolved by the explicit form of $\tau_k$, which does not contain any $\rho$ dependence in the $B$-direction near $\rho=0$. 
This extended map is still Hamiltonian, and satisfies the twist condition for $k\leq \ell$, but not for $k>\ell$.

Therefore it satisfies the claim of the proposition.
\end{proof}

\begin{remark}
Given a return map $\tau$ that is Hamiltonian, we point out that the Hamiltonian family generating $\tau$ is not unique, and more importantly, that various dynamical properties depend on the choice of Hamiltonian. For example, on the disk $(D^2, rdr\wedge d\theta)$, the return map $\tau =\id$ is generated by the autonomous Hamiltonians $H_k=k\pi r^2$. For given $k$, the Robbin-Salamon index of the $1$-periodic orbit at $\partial D^2$ is $2k$, i.e.\ $k$-dependent. The associated paths of symplectic matrices have the same endpoints, but are not homotopic rel endpoints. This also illustrates the interpretation of the RS-index as a winding number.
Note that $D^2$ has a Hamiltonian circle action that extends over the whole space. We don't know whether the same type of phenomenon occurs for more general symplectic manifolds (i.e.\ without a global Hamiltonian circle action).
\end{remark}

\subsection{Non-examples: Katok examples}\label{app:Katok}
In \cite{K73}, Katok constructed examples of non-reversible Finsler metrics on $S^n$ with only finitely many simple closed geodesics. Here is a description of such examples using Brieskorn manifolds.
We consider 
$$
\Sigma^{2n-1} := \left\{ (z_0,\ldots,z_n) \in \C^{n+1} ~\Bigg|~\sum_j z_j^2 = 0 \right\}
\cap 
S^{2n+1}_1,
$$
equipped with the contact form $\alpha =\frac{i}{2}\sum_j z_j d\bar z_j -\bar z_j dz_j$.
These spaces are contactomorphic to $S^*S^n$ with its canonical contact structure. The given contact form is actually the prequantization form.

We describe the setup in detail when $n=2m+1$ is odd.
We group the coordinates in pairs, and make the following unitary coordinate transformation:
$$
w_0=z_0,
w_1=z_1,
w_{2j}=\frac{\sqrt 2}{2} (z_{2j}+i z_{2j+1}),
w_{2j+1}=\frac{i \sqrt 2}{2} (z_{2j}-i z_{2j+1})
\text{ for }j=1,\dots,m.
$$
Because this is a unitary transformation, the form $\alpha$, expressed in $w$-coordinates, still has the form
$$
\alpha=\frac{i}{2} \sum_j w_j d\bar w_j -\bar w_j dw_j
$$
For a tuple $\epsilon=(\epsilon_1,\dots, \epsilon_m) \in (-1,1)^m$, define the function $H_\epsilon$ on a neighborhood of $\Sigma^{2n-1}$ via
$$
H_\epsilon(w)=\Vert w \Vert^2+\sum_{j} \epsilon_j( 
|w_{2j}|^2-|w_{2j+1}|^2
).
$$
For $\epsilon$ sufficiently small, this function is positive, so we define a perturbed contact form by
$$
\alpha_\epsilon = H_\epsilon^{-1} \cdot \alpha.
$$
The Reeb vector field of $\alpha_\epsilon$ is 
$$
R_\epsilon= X_\epsilon+\overline{X}_\epsilon,$$
where
$$
X_\epsilon=i w_0\frac{\partial}{\partial w_0} + i w_1\frac{\partial}{\partial w_1}
+\sum_j\left( i(1+\epsilon_j) \frac{\partial }{\partial w_{2j}} +  i(1-\epsilon_j) \frac{\partial }{\partial w_{2j+1}}\right),
$$
$$
\overline{X}_\epsilon=-i \overline{w}_0\frac{\partial}{\partial \overline{w}_0} - i \overline{w}_1\frac{\partial}{\partial \overline{w}_1}
-\sum_j\left( i(1+\epsilon_j) \frac{\partial }{\partial \overline{w}_{2j}} +  i(1-\epsilon_j) \frac{\partial }{\partial \overline{w}_{2j+1}}\right).
$$
The Reeb flow is therefore given by
$$
(w_0,\dots,w_n)\longmapsto
(e^{2\pi it}w_0,e^{2\pi it}w_1,e^{2\pi it(1+\epsilon_1)}w_2,e^{2\pi it(1-\epsilon_1)}w_3,\dots, e^{2\pi it(1+\epsilon_m)}w_{n-1},e^{2\pi it(1-\epsilon_m)}w_{n}).
$$
This flow has only $n+1$ periodic orbits if all $\epsilon_j$ are rationally independent. These are given by
\[
\begin{split}
\gamma_0(t)=& \left(\frac{1}{\sqrt{2}}e^{2\pi it},\frac{i}{\sqrt{2}}e^{2\pi it},0,\dots,0\right), t \in [0,1]\\
\beta_0(t)=&\left(\frac{1}{\sqrt{2}}e^{2\pi it},-\frac{i}{\sqrt{2}}e^{2\pi it},0, \dots, 0\right), t \in [0,1]\\
\gamma_j(t)=&\left(0,0,\dots, e^{2\pi it(1+\epsilon_j)},0,\ldots,0,0\right), t \in [0,1/(1+\epsilon_j)]\\
\beta_j(t)=&\left(0,0,\dots, 0, e^{2\pi it(1-\epsilon_j)},\ldots,0,0\right), t \in [0,1/(1-\epsilon_j)]\\
\end{split}
\]
for $j=1,\dots,m$. 
\begin{remark}
As stated, we see that there are only finitely many periodic orbits.
Furthermore, since the unperturbed system, i.e.~$\epsilon=0$, describes the geodesic flow on the round sphere, and the perturbation $\alpha_\epsilon$ is $C^2$-small for small $\epsilon$, it follows that the Reeb flow of the contact form $\alpha_\epsilon$ corresponds to the geodesic flow of a Finsler metric.
In Section~\ref{sec:katok_examples_relation} we describe how to obtain an explicit relation with the famous Katok examples for $S^*S^2$.
\end{remark}
We construct a supporting open book for the contact form $\alpha_\epsilon$ using the map
$$
\Theta: \Sigma^{2n-1} \longrightarrow \C, (w_0,w_1,\ldots,w_n) \longmapsto w_0.
$$
The zero set of $\Theta$ defines the binding, the pages are the sets of the form $P_\theta=\{\arg \Theta=\theta\}$, $\theta \in S^1$, which are all copies of $\mathbb{D}^*S^{n-1}$, and the monodromy is $\tau^2$ where $\tau$ is the Dehn-Seidel twist. The (boundary extended) return map for the page $P_0=\Theta^{-1}(\R_{>0})\cong \mathbb{D}^*S^{n-1}$ is
\[
\begin{split}
\Phi: P_0&\longrightarrow P_0, \\
p=(r_0,w_1,w_2,w_3,\ldots,w_{n-1},w_{n}) & \longmapsto 
(r_0, w_1,e^{2\pi i\epsilon_1}w_2,e^{-2\pi i \epsilon_1}w_3,
\ldots,\\
&\phantom{\longmapsto (~} e^{2\pi i\epsilon_{m}}w_{n-1}, e^{-2\pi i\epsilon_{m}}w_n).
\end{split}
\]
Here, $w_0=r_0 \in \mathbb{R}_{\geq 0}$ is a real non-negative number, and note that the first return time is constant equal to $1$ (which follows by looking at the first coordinate). If all $\epsilon_j$ are irrational and rationally independent, this map has only two periodic points, both actually fixed, given by
\[
\begin{split}
p_0=& \left(\frac{1}{\sqrt{2}},\frac{i}{\sqrt{2}},0,\dots,0\right)\\
q_0=&\left(\frac{1}{\sqrt{2}},-\frac{i}{\sqrt{2}},0, \dots, 0\right).\\
\end{split}
\]
Note that $p_0,q_0$ are both interior fixed points, and irrationality of the $\epsilon_j$ implies that there are no boundary fixed points.
We will explain now why this map is Hamiltonian with boundary preserving Hamiltonian flow.
The symplectic form on the interior of the page $P_0$ is the restriction of $d\alpha_{\epsilon}$.
To manipulate this, let us define
$$
H= \Vert w \Vert^2,\quad \Delta_\epsilon =\sum_{j} \epsilon_j( 
|w_{2j}|^2-|w_{2j+1}|^2),
$$
so $H_\epsilon=H+\Delta_\epsilon$.
Observe that the return map $\Phi$ is generated by the $2\pi$-flow of the vector field
$$
X=i\sum_{j=1}^m\epsilon_j\left(
w_{2j}\frac{\partial}{\partial w_{2j}}-\bar w_{2j}\frac{\partial}{\partial \bar w_{2j}}
-
w_{2j+1}\frac{\partial}{\partial w_{2j+1}}+\bar w_{2j+1}\frac{\partial}{\partial \bar w_{2j+1}}
\right).
$$
This vector field is tangent to the page and preserves $H$ and $\Delta$, and hence also $H_\epsilon$.
Plug $X$ in into $d\alpha_\epsilon$. We find
\[
\begin{split}
\iota_X(d H_\epsilon^{-1}\wedge \alpha +H_\epsilon^{-1} d\alpha )& =
-\alpha(X) d H_\epsilon^{-1}+H_{\epsilon}^{-1} \iota_X d\alpha\\
&=-\Delta_\epsilon d H_\epsilon^{-1} -H_\epsilon^{-1} d \Delta_\epsilon= -d (H_\epsilon^{-1} \Delta_\epsilon).
\end{split}
\]
This means that the Hamiltonian generating the return map is $H_\epsilon^{-1} \Delta_\epsilon$.
Moreover, index-positivity follows, by observing that it holds for the round metric on $S^2$ and the fact that it is an open condition. It follows from Thm.\ \ref{thmB} that $\Phi$ does \emph{not} satisfy the twist condition for any Liouville structure on $\mathbb{D}^*S^2$.

\begin{remark}\label{rk:evenn}
The setup for $n$ even is very similar: we drop the $w_0$-coordinate.
\end{remark}

\subsection{Relation with the Katok examples}
\label{sec:katok_examples_relation}
We explain how to see that the above dynamical systems indeed correspond to the Katok examples in case of $S^*S^2$ (i.e.\ $n=2$).
More precisely, we will show that the geodesic flow of the Katok examples is conjugated to the Reeb flow of $\alpha_\epsilon$.
We need some preparation, which applies to all dimensions, before we specialize to dimension 3.
First of all, we fix positive weights $(a_1,\ldots, a_n)\in \R_{>0}^n$.
Then we define the $1$-forms on the sphere $S^{2n-1}$ given by
$$
\beta_0
=\iota^*\left(\frac{\sum_j x_j dy_j-y_j dx_j}{\sum_k a_k (x_k^2+y_k^2)}\right)=\iota^*\left(\frac{1}{\sum_k a_k \vert z_j\vert^2}\frac{i}{2}\sum_j z_j d\bar z_j-\bar z_j dz_j\right)
$$
and
$$
\beta_1=
\iota^*\left(\sum_j \frac{1}{a_j}(x_j dy_j-y_j dx_j)\right)=\iota^*\left(\frac{i}{2}\sum_j \frac{1}{a_j}\left(z_j d\bar z_j-\bar z_j dz_j\right)\right),
$$
where $\iota$ is the inclusion map $S^{2n-1}\to \R^{2n}$.
We will show that the first form is a contact form and that it is strictly contactomorphic to the latter. For this, consider the map
\[
\begin{split}
\psi: S^{2n-1} & \longrightarrow S^{2n-1} \\
(z_1,\ldots,z_n) & \longmapsto
\left(
\sqrt{ \frac{a_1}{\sum_k a_k(x_k^2+y_k^2)} }z_1,\ldots,
\sqrt{ \frac{a_n}{\sum_k a_k(x_k^2+y_k^2)} }z_n
\right).
\end{split}
\]
We find
\[
\begin{split}
\psi^*\beta_1
&=
\sum_j \frac{1}{a_j}\sqrt{ \frac{a_j}{\sum_k a_k(x_k^2+y_k^2)}}\left(\sqrt{ \frac{a_j}{\sum_k a_k(x_k^2+y_k^2)}}
(
x_jdy_j-y_jdx_j
)
+(x_jy_j -y_jx_j) d\left( \sqrt{ \frac{a_j}{\sum_k a_k(x_k^2+y_k^2)} } \right)\right) \\
&=\beta_0.
\end{split}
\]
This also shows that $\beta_0$ is a contact form, as $\psi$ is a diffeomorphism. We have shown:

\begin{lemma}
\label{lemma:weighted_form}
The form $\beta_0$ is a contact form, and it is strictly contactomorphic to $\beta_1$.
The Reeb field for $\beta_k$ for $k=0,1$ is given by
$$
R=\sum_j a_j \left(x_j \frac{\partial}{\partial y_j}-y_j \frac{\partial}{\partial x_j}
\right)
.
$$
\end{lemma}

We now specialize to the $3$-dimensional situation, for which in $w$-coordinates (cf.\ Remark \ref{rk:evenn}) we have
$$
\Sigma^3=\left\{(w_1,w_2,w_3)\in \mathbb{C}^3:w_1^2-2iw_2w_3=0\right\}\cap S^3.
$$
Consider the explicit covering map
\[
\pi: S^3 \longrightarrow \Sigma^3,
(z_0,z_1) \longmapsto \left(w_1=\sqrt{2}z_0 z_1,w_2=z_0^2,w_3=-iz_1^2\right).
\]
We quickly verify that this is a covering map:
\begin{itemize}
\item we have $w_1^2-2i w_2w_3=2z_0^2z_1^2-2z_0^2z_1^2=0$;
\item we have $|w_1|^2+|w_2|^2+|w_3|^2=2|z_0|^2|z_1|^2+|z_0|^4+|z_1|^4=(|z_0|^2+|z_1|^2)^2=1$;
\item the map is two to one, since all entries are quadratic. 
\end{itemize}
We compute the pullback $\pi^*\alpha_\epsilon$,
\[
\begin{split}
\pi^*\alpha_\epsilon
&=
\frac{1}{2|z_0|^2|z_1|^2+|z_0|^4+|z_1|^4 +\epsilon |z_0|^4-\epsilon |z_1|^4}
2(|z_0|^2+|z_1|^2) \frac{i}{2}\sum_j (z_j d\bar z_j -\bar z_j d z_j)\\
&=
\frac{2(|z_0|^2+|z_1|^2)}{(|z_0|^2+|z_1|^2)((1+\epsilon) |z_0|^2 +(1-\epsilon) |z_1|^2)}
\frac{i}{2}\sum_j (z_j d\bar z_j -\bar z_j d z_j)\\
&= \frac{2}{((1+\epsilon) |z_0|^2 +(1-\epsilon) |z_1|^2)}
\frac{i}{2}\sum_j (z_j d\bar z_j -\bar z_j d z_j).
\end{split}
\]
By Lemma~\ref{lemma:weighted_form}, the form $\pi^* \alpha_\epsilon$ to strictly contactomorphic to the contact form $\beta_1$ with weights $a_1=1+\epsilon$ and $a_2=1-\epsilon$, which is just the ellipsoid model for the contact $3$-sphere.
To complete the argument, we use a result due to Harris and Paternain \cite[Section 5]{HP}, which relates the ellipsoids to the Katok examples.

\section{Symplectic homology of surfaces}
\label{app:SH_surfaces}
Let us consider connected Liouville domains in dimension $2$.
The simplest such Liouville domain is $D^2$, which has vanishing symplectic homology.
For all other surfaces, note:
\begin{lemma}
Let $(W,\lambda)$ be a connected Liouville domain of dimension $2$. Assume that $W$ is not diffeomorphic to $D^2$. Take a periodic Reeb orbit $\delta$ on one of the boundary components of $W$.
Then $[\delta] \in \tilde \pi_1(W)$ is non-trivial.
Furthermore, if $\delta_1$ and $\delta_2$ are periodic Reeb orbits on different boundary components, then $[\delta_1]\neq [\delta_2]$ as free homotopy classes. $\hfill \square$
\end{lemma}

Assume $W\neq D^2$, and denote the completion by $\widehat W$. Then the chain complex for an admissible Hamiltonian $\widehat H$ that is both negative and $C^2$-small on $W$ has the form
$$
CF_\bullet(\widehat H)
=\bigoplus_{\delta \in \tilde \pi_1(W)} CF^{\delta}_\bullet(\widehat H),
$$
where $CF_\bullet^{\delta}(\widehat H)$ is generated by $1$-periodic orbits in the free homotopy class $\delta$. The direct summand corresponding to contractible orbits needs as least as many generators as $\rk H_\bullet(W)$ by the Morse inequalities.
\\
\noindent

\begin{lemma}
For each class $\delta$ the direct summand $CF_\bullet^{\delta}(\widehat H)$ forms a subcomplex, and so we have a splitting
$$
HF_\bullet(\widehat H)=\bigoplus_{\delta \in \tilde \pi_1(W)} HF_\bullet^\delta(\widehat H).
$$ 
In addition, as ungraded modules we have
\[
HF_\bullet^\delta(\widehat H)
\cong \begin{cases}
\Z^2 & \text{if }\delta \text{ is a positive boundary class}, \text{ and }\mbox{slope}(\widehat H)  \text{ is sufficiently large},\\
H_\bullet(W) & \text{if }\delta \text{ is the trivial class,}\\
0 & \text{otherwise.}
\end{cases}
\]
\end{lemma}

Here, a positive boundary class just means a homotopy class of a positive multiple of a boundary component (oriented according the positive boundary orientation).

\begin{proof}
The first assertion follows from the fact that Floer cylinders do not change the free homotopy class.
For the second claim we use:
\begin{itemize}
\item The Floer differential of a $C^2$-small Hamiltonian between critical points is the Morse differential, which implies the second case. 
\item After a suitable Morse perturbation (a Morse function on $S^1$ with precisely two critical points) breaking the $S^1$-symmetry given by time-shifts, each positive boundary class gives two generators, corresponding to the critical points of the Morse function on $S^1$; as shown in \cite{CFHW}, the differential is the Morse differential, which vanishes. Moreover, this symmetry-breaking process preserves the homotopy classes of periodic orbits, as observed in Remark~\ref{rem:free_homotopy}.
\end{itemize}
\end{proof}

\begin{corollary}
\label{cor:SH_surface}
Suppose that $W$ is a connected Liouville domain of dimension $2$. Assume that $W$ is not diffeomorphic to $D^2$.
Then as an ungraded module we have
\[
SH_\bullet(W)
\cong
H_\bullet(W)
\oplus
\bigoplus_{\delta \text{ positive boundary class}} \Z^2.
\]
\end{corollary} 
\section{On symplectic return maps}\label{app:longorbits}

In this appendix, for convenience of the reader, we collect some standard facts concerning return maps arising from a given Reeb dynamics on some contact manifold (cf.\ the construction of the \emph{Calabi homomorphism}, e.g.\ in \cite[Sec.\ 10.3]{MS17}, or \cite[Sec.\ 3.3]{ABHS} for the case of the $2$-disk). In particular, we show that long Hamiltonian orbits on a global hypersurface of section correspond to long Reeb orbits on the ambient contact manifold.

Consider a map $\tau:\mbox{int}(\Sigma) \rightarrow \mbox{int}(\Sigma)$ defined on the interior of a $2n$-dimensional Liouville domain $\Sigma$. We assume that $\Sigma$ arises as a (connected) global hypersurface of section for some Reeb dynamics on a $2n+1$-dimensional contact manifold $(M,\alpha)$, and $\tau$ is the associated return map. Let $R_\alpha$ be the Reeb vector field of $\alpha$. Denote by $B=\partial \Sigma$, which we assume to be a contact submanifold of $M$ with induced contact form $\alpha_B=\alpha\vert_B$, so that $R_\alpha\vert_B$ is tangent to $B$. Let $\lambda=\alpha\vert_\Sigma$, which is a Liouville form on $\mbox{int}(\Sigma)$ since $R_\alpha$ is assumed to be positively transverse to the interior of $\Sigma$. That is, the two-form $\omega=d\lambda$ is symplectic on $\mbox{int}(\Sigma)$. The $1$-form $\lambda_B=\lambda\vert_B$ coincides with the contact form $\alpha_B$. Note that it is degenerate along $B$. By Stokes' theorem, the symplectic volume of $\Sigma$ then coincides with the contact volume of $B$:
$$
\mbox{vol}(\Sigma,\omega)=\int_\Sigma \omega^n=\int_\Sigma d(\lambda \wedge d\lambda^{n-1})=\int_B \alpha_B \wedge d\alpha_B^{n-1}=\mbox{vol}(B,\alpha_B).
$$

Note that $\tau$ is automatically a symplectomorphism with respect to $\omega$. Indeed, denote the time-$t$ Reeb flow by $\varphi_t$, and let $T:\mbox{int}(\Sigma)\rightarrow \mathbb{R}^+$
$$
T(x)=\min\{t>0:\varphi_t(x)\in \mbox{int}(\Sigma)\}
$$
denote the first return time function. Then $\tau(x)=\varphi_{T(x)}(x)$, and so, for $x \in \mbox{int}(\Sigma)$, $v \in T_x\Sigma$, we have
$$
d_x\tau(v)=d_xT(v)R_{\alpha}(\tau(x))+d_x\varphi_{T(x)}(v).
$$
Using that $\varphi_t$ satisfies $\varphi_t^*\alpha=\alpha$, we obtain
\begin{equation}
    \begin{split}
        (\tau^*\lambda)_x(v)&=\alpha_{\tau(x)}(d_x\tau(v))\\
        &=d_xT(v)+(\varphi_{T(x)}^*\alpha)_x(v)\\
        &=d_xT(v)+\lambda_x(v).\\
    \end{split}
\end{equation}
Therefore 
\begin{equation}\label{Tlambda}
\tau^*\lambda=dT + \lambda,
\end{equation}
which in particular implies that $\tau^*\omega=\omega$. 

Moreover, the average of the return time function gives the contact volume of $M$, i.e.\ we have the identity
\begin{equation}\label{contactvol}
\int_\Sigma T \omega^n=\mbox{vol}(M,\alpha).
\end{equation}
This may be proved as follows. We have a smooth embedding
$$
\psi \colon \mathbb{R}/\mathbb{Z}\times \mbox{int}(\Sigma) \rightarrow M,
$$
given by $\psi(s,x)=\varphi_{sT(x)}(x)$, which is a diffeomorphism onto $M\backslash B$. It satisfies $$(\psi^*\alpha)(\partial_s)=\alpha(T R_\alpha)=T,$$
and, for $v \in T\mbox{int}(\Sigma)$, 
$$(\psi^*\alpha)(v)=\alpha(sdT(v)R_\alpha +d\varphi_{sT}(v))=sdT(v)+\alpha(v).$$
Then
$$
\psi^*\alpha=T ds + sdT+\lambda=d(sT)+\lambda,
$$
and so
$$
\psi^*(\alpha \wedge d\alpha^n)=(d(sT)+\lambda)\wedge d\lambda^n=T ds \wedge \omega^n.
$$
Integrating, and using the fact that $B$ is codimension 2 in $M$, we obtain
$$
\mbox{vol}(M,\alpha)=\int_{M\backslash B}\alpha \wedge d\alpha^n=\int_{\mathbb{R}/\mathbb{Z}\times \mbox{int}(\Sigma)}\psi^*(\alpha \wedge d\alpha^n)
$$
$$
=\int_{\mathbb{R}/\mathbb{Z}\times \mbox{int}(\Sigma)}T ds \wedge \omega^n=\int_{\mbox{int}(\Sigma)}T \omega^n=\int_{\Sigma}T \omega^n,
$$
where we have used that $\omega^n\vert_B\equiv 0$, and the claim follows. In case where $\tau$ is Hamiltonian, we want to relate the Hamiltonian action of a periodic orbit of $\tau$ to the Reeb action of the corresponding Reeb orbit in the ambient contact manifold.

Let $H: S^1 \times \Sigma \rightarrow \mathbb{R}^+$ be a Hamiltonian generating $\tau$, i.e.\ the isotopy $\phi_t$ defined by $\phi_0=id$, $\frac{d}{dt}\phi_t=X_{H_t}\circ \phi_t$ satisfies $\phi_1=\tau$. The sign convention for the Hamiltonian vector field is $i_{X_{H_t}}\omega=-dH_t$. We usually view this Hamiltonian isotopy as defining an element $\phi=\phi_H=[\{\phi_t\}]$ in the universal cover $\widetilde{\mbox{Diff}}(\Sigma,\omega)$ of the space of symplectomorphisms $\mbox{Diff}(\Sigma,\omega)$. By Cartan's formula, we have
$$
\partial_t \phi_t^*\lambda=\phi_t^*\mathcal{L}_{X_{H_t}}\lambda=\phi_t^*(i_{X_{H_t}}\omega + d(i_{X_{H_t}}\lambda))=\phi_t^*d(i_{X_{H_t}}\lambda-H_t),
$$
and so integrating we obtain
\begin{equation}\label{Tlambda2}
\tau^*\lambda-\lambda=dF_H,
\end{equation}
where 
\begin{equation}\label{FHdef}
F_H=\int_0^1 (i_{X_{H_t}}\lambda - H_t) \circ \phi_t \;dt
\end{equation}
Combining (\ref{Tlambda}) and (\ref{Tlambda2}) we deduce that 
\begin{equation}\label{Ftau}
\tau=F_H+C
\end{equation}
for some constant $C$ (assuming $\Sigma$ is connected). 

We determine the constant $C$ under a suitable assumption, which we assume holds in all what follows. Namely, assume that $\tau$ extends to $\Sigma$ with the same formula, i.e.\ via an extension of the return time function $T$ to $\Sigma$. Assume also that $H_t\vert_B\equiv const:=C_t>0$ for some $H$ generating $\tau$. Equivalently, $X_{H_t}\vert_B=h_tR_B$ for some (not necessarily positive) smooth function $h_t$ on $B$, satisfying $h_t=dH_t(V_\lambda)\vert_B$ where $V_\lambda$ is the Liouville vector field associated to $\lambda$. In this case, denoting $\gamma_x(t)=\phi_t(x)$ for $x \in B$ and $t\in [0,1]$, we get
\begin{equation}\label{FH}
F_H(x)=\int_{\gamma_x} \lambda_B-\int_0^1C_tdt=\int_0^1 (h_t(\phi_t(x))-C_t)dt,
\end{equation}
On the other hand, let $\beta_x(t)=\varphi_t(x)$ be the Reeb orbit through $x$ ending at $\beta_x(1)=\tau(x)$, for $t \in [0,1]$, which we assume parametrized so that $\dot{\beta}_x=T(x)R_B(\beta_x)$. Note that $\beta_x$ is a reparametrization of $\gamma_x$, and so we obtain
$$
\tau(x)=\int_{\beta_x} \lambda_B=\int_{\gamma_x}\lambda_B
$$
This means that $T$ is the unique primitive of $\tau^*\lambda-\lambda$ satisfying $T(x)=\int_{\gamma_x}\lambda_B$ for $x \in B$. Combining (\ref{Ftau}) and (\ref{FH}), we conclude that
$$
C=\int_0^1C_tdt>0,
$$
a positive constant. 

By the above computation, $T$ is what is usually called the action of $\phi=\phi_H$ with respect to $\lambda$, and is independent of the isotopy class (with fixed endpoints) of the path $\phi_H$. The Calabi invariant is then by definition the average action $CAL(\phi_H,\omega)=\int T \omega^n$, which is independent of $\lambda$; cf.\ \cite{MS17,ABHS}. Combining with (\ref{contactvol}), we obtain
$$
CAL(\phi_H,\omega)=\mbox{vol}(M,\alpha).
$$

Let $\gamma:S^1=\mathbb{R}/k\mathbb{Z}\rightarrow \Sigma$, defined by $\gamma(t)=\phi_t(x)$, be a $k$-periodic Hamiltonian orbit associated to the $k$-periodic point $x$ of $\tau$. That is, we have $x=\gamma(0)$, $\gamma(1)=\tau(x),\dots, \gamma(k)=\tau^k(x)=x$, and assume that $k$ is the minimal period of $x$. We then get
$$
\sum_{i=1}^k F_H(\tau^i(x))=\mathcal{A}_{H^{\#k}}(\gamma)
$$
is precisely the Hamiltonian action of $\gamma$ with respect to the Hamiltonian 
        $$
        H_t^{\#k}=\sum_{i=1}^k H_t \circ \phi_t^{-i}
        $$
generating $\tau^k$. If $\beta: S^1=\mathbb{R}/\mathbb{Z} \rightarrow M$ is the Reeb orbit corresponding to $\gamma$, (\ref{Ftau}) implies that its period is
$$
\int_{S^1}\beta^*\alpha=\sum_{i=1}^k T(\tau^i(x))= \sum_{i=1}^k F_H(\tau^i(x))+kC=\mathcal{A}_{H^{\#k}}(\gamma)+kC
$$
Since $C>0$, this implies the following: if the Hamiltonian action of every $k$-periodic orbit $\gamma$ grows to infinity with $k$, then the period of the associated Reeb orbits $\beta$ also. In other words, long Hamiltonian periodic orbits in the global hypersurface of section give long Reeb orbits in the ambient contact manifold. 

We summarize the above discussion in the following:

\begin{lemma}\label{lemma:appendix} Let $(M^{2n+1},\alpha)$ be a contact manifold, $(\Sigma^{2n},\omega=d\alpha\vert_\Sigma)$ a Liouville domain which is a global hypersurface of section for the Reeb flow, $(B^{2n-1},\alpha_B)=(\partial \Sigma,\alpha\vert_B)$, $\tau:\mbox{int}(\Sigma) \rightarrow \mbox{int}(\Sigma)$ the Poincar\'e return map, and $T: \mbox{int}(\Sigma)\rightarrow \mathbb R^+$ the first return time. Then:
\begin{enumerate}
    \item $\mbox{vol}(\Sigma,\omega)=\mbox{vol}(B,\alpha_B)$.
    \item $\mbox{vol}(M,\alpha)=\int_{\Sigma}T \omega^n$.
    \item $\tau$ is an exact symplectomorphism.
    \item If $\tau$ is Hamiltonian with generating isotopy $\phi_H=[\{\phi_t\}]\in \widetilde{\mbox{Diff}}(\Sigma,\omega)$, and extends to $\Sigma$ as a (not necessarily positive) reparametrization of the Reeb flow at $B$, then:
    \begin{itemize}
        \item[(i)] $CAL(\phi_H,\omega)=\mbox{vol}(M,\alpha)$.
        \item[(ii)] The period of a Reeb orbit $\beta$ on $M$ corresponding to a $k$-periodic Hamiltonian orbit $\gamma$ on $\Sigma$ is
        $$\int_{S^1}\beta^*\alpha=\mathcal{A}_{H^{\#k}}(\gamma)+kC$$ for some positive constant $C>0$, where 
        $$
        \mathcal{A}_{H^{\#k}}(\gamma)=\int_{S^1}\gamma^*\lambda-\int_0^1H^{\#k}_t(\gamma(t))dt
        $$
        is the Hamiltonian action of $\gamma$ with respect to the Hamiltonian 
        $$
        H_t^{\#k}=\sum_{i=1}^k H_t \circ \phi_t^{-i}
        $$
        generating $\tau^k$. In particular, if $\gamma$ has large action, then $\beta$ has large period.
    \end{itemize}
\end{enumerate}
\end{lemma}

\section{Strong convexity implies strong index-positivity}\label{app:indexpos}

In this appendix, we give a general condition for index-positivity to hold, which is also relevant for the restricted three-body problem. A connected compact hypersurface $\Sigma \subset \mathbb{R}^{4}$ is said to \emph{bound a
strongly convex domain} $W \subset \mathbb{R}^{4}$ whenever there exists a smooth function $\phi: \mathbb{R}^{4}\rightarrow \mathbb{R}$ satisfying:
\begin{itemize}
    \item[(i)](Regularity) $\Sigma=\{\phi=0\}$ is a regular level set;
    \item[(ii)](Bounded domain) $W = \{z\in \mathbb{R}^{4}: \phi(z) \leq 0\}$ is bounded and contains the origin; and
    \item[(iii)](Positive-definite Hessian) $\nabla^2\phi_{z}(h, h) > 0$ for $z\in W$ and for each non-zero tangent vector $h\in T\Sigma$.
\end{itemize}
In this case, the radial vector field is transverse to $\Sigma$, and so $\Sigma$ is a contact-type $3$-sphere, inheriting a contact form $\alpha$ induced by the standard Liouville form in $\mathbb{R}^{4}$. 

\begin{lemma}\label{lemma:indexposapp}
Suppose that $\Sigma$ bounds a strongly convex domain. Then $\Sigma$ is strongly index-positive.
\end{lemma}

\begin{remark}
In the planar restricted three-body problem, the values of energy/mass ratio $(c,\mu)$ for which the Levi-Civita regularization bounds a strictly convex domain is called the \emph{convexity range}, which in particular implies that the dynamics is \emph{dynamically convex} (cf.\ \cite{HWZ98,AFFHvK, AFFvK}). It follows that index-positivity holds in the convexity range for the quotient $\mathbb{R}P^3$, which is part of the assumptions of Thm.\ \ref{thmA}.
\end{remark}

\begin{proof}
Write $\Sigma= \phi^{-1}(0)$ as in the definition above. Denote the contact form on $\Sigma$ by $\alpha:=\lambda|_{\Sigma}$.
We will use the standard quaternions $I,J,K$, where $I$ is chosen to coincide with the standard complex structure.

The tangent space of $\Sigma$ is spanned by the vectors 
$$
R=X_\phi /\alpha(X_\phi) =I \nabla \phi /\alpha(X_\phi) =Iw,\; 
U=Jw -\alpha(Jv)R,\;
V=Kw -\alpha(Kv)R.
$$
We note that $U$ and $V$ give a symplectic trivialization $\epsilon$ of $(\xi=\ker \alpha,d\alpha)$.
To see this, we compute
$$
d\alpha(U,V)= d \alpha ( Jw, Kw) =w^t J^t I^t K w=w^t K^t K w=w^t w=1.
$$
In order to prove the claim, we investigate the rate of change of a version of the rotation number. 
See Ch.\ 10.6 in \cite{FvK18} for a detailed standard description of the Robbin-Salamon index in terms of the rotation number.
We will detail the version that we will use below.

We look at the linearization of the Hamiltonian flow:
\begin{equation}
\label{eq:linearized_flow_cvx_vf}
\dot X =\nabla_X X_\phi= I \nabla^2\phi \cdot X.
\end{equation}
Starting with $X(0) \in \xi$, we compute how quickly the vector $X$ rotates with respect to the frame. Define the angular form
$$
\Theta=
\frac{udv-vdu}{u^2+v^2}
=
\frac{u d\alpha(U, \cdot) +v d\alpha(V,\cdot) }{u^2+v^2}
=\frac{d\alpha(u U+v V,\cdot)}{u^2+v^2}
,
$$
where $(u,v)$ are cartesian coordinates on the plane spanned by the frame $(U,V)$, so we may write $X=u U+ v V$.
We plug in $\dot X$ and find
\begin{equation}
\label{eq:angle_change}
\Theta(\dot X)=\frac{d\alpha(X,\dot X)}{u^2+v^2}=
\frac{(u U+v V)^t I^t I \nabla^2\phi \cdot (u U+v V)}{u^2+v^2}
=
\frac{\nabla^2\phi(uU+vV,uU+vV)}{u^2+v^2} \geq \lambda_{min}>0,
\end{equation}
where $ \lambda_{min}$ is the minimal eigenvalue of $\nabla^2\phi$ over the compact hypersurface $\Sigma$.
After we have set up some notation, we will see that this is enough to get a lower bound on the growth rate of the Robbin-Salamon index.
With our global trivialization $\epsilon$, we can define the matrix
$$
\psi(t)=\epsilon \circ dFl^R_t \circ \epsilon^{-1}.
$$
By applying Equation~\ref{eq:linearized_flow_cvx_vf} to the initial vectors $\epsilon^{-1}(1,0)$ and $\epsilon^{-1}(0,1)$, we get a linear evolution equation for the matrix $\psi(t)$,
\begin{equation}
\label{eq:psi_evolution}
\dot \psi =A(t) \psi,
\end{equation}
where $A$ is a time-dependent matrix.
We will view this ODE as a vector field on $Sp(2)$: the linearized Reeb flow along each Reeb orbit will give rise to such a vector field.

To relate the above angle to the Conley-Zehnder index, we also need to recall the Iwasawa decomposition, also known as KAN decomposition, of $Sp(2)$. 
Write 
$$
KAN:=\left\{\left( 
\left(
\begin{array}{cc}
\cos(\phi) & -\sin(\phi) \\
\sin(\phi) & \cos(\phi)
\end{array}
\right),
\left(
\begin{array}{cc}
a & 0 \\
0 & a^{-1}
\end{array}
\right)
,
\left(
\begin{array}{cc}
1 & t \\
0 & 1
\end{array}
\right)\right)
~\Bigg|~
\phi \in [0,2\pi), a \in \R_{>0},t \in \R
\right\}
$$
And put
\[
\mbox{kan}: 
KAN \longrightarrow Sp(2),
(\phi,a, t) \longmapsto 
\left(
\begin{array}{cc}
\cos(\phi) & -\sin(\phi) \\
\sin(\phi) & \cos(\phi)
\end{array}
\right)
\left(
\begin{array}{cc}
a & 0 \\
0 & a^{-1}
\end{array}
\right)
\left(
\begin{array}{cc}
1 & t \\
0 & 1
\end{array}
\right).
\]
This map has the inverse
\[
\begin{split}
\mbox{kan}^{-1}: Sp(2) &\longrightarrow KAN \\
\left(
\begin{array}{cc}
a & b \\
c & d
\end{array}
\right) &\longmapsto\left(
\frac{1}{\sqrt{a^2+c^2}} 
\left(
\begin{array}{cc}
a & -c \\
c & a
\end{array}
\right)
,
\left(
\begin{array}{cc}
\sqrt{a^2+c^2} & 0 \\
0 & \frac{1}{\sqrt{a^2+c^2}} 
\end{array}
\right)
,
\left(
\begin{array}{cc}
1 & \frac{ab+cd}{\sqrt{a^2+c^2}} \\
0 & 1
\end{array}
\right)\right).
\end{split}
\]
The KAN angle can locally be determined as
$$
\arg( \mbox{kan}^{-1}(\psi)) =\atan(c/a),
$$
so we see that the change in angle equals
$$
\frac{d}{dt}\arg( \mbox{kan}^{-1}(\psi)) =\frac{d}{dt}\atan(c/a)=\frac{a\dot c -c\dot a}{a^2+c^2}.
$$
On the other hand, the rate of change of the KAN angle equals $\Theta(\dot X)$. 
Indeed, the first column of $\psi(t)$ is the vector $Z(t):=\epsilon(X(t))=\left( \begin{array}{c} u \\ v \end{array} \right)$ if we put $X(0)=\epsilon^{-1}(1,0)$.

By Equation~\eqref{eq:angle_change}, this rate of change is at least $\lambda_{min}$, where $\lambda_{min}$ is the minimal eigenvalue of $-IA$ (which we assume to be positive definite). This means that each slice
$$
S_\phi=\mbox{kan}\left(\left\{ 
\left(\left(
\begin{array}{cc}
\cos(\phi) & -\sin(\phi) \\
\sin(\phi) & \cos(\phi)
\end{array}
\right),
\left(
\begin{array}{cc}
a & 0 \\
0 & a^{-1}
\end{array}
\right)
,
\left(
\begin{array}{cc}
1 & t \\
0 & 1
\end{array}
\right)\right)
~\Bigg|~
 a \in \R_{>0},t \in \R
\right\}\right)
$$
is a global surface of section for the vector field associated with Equation~\eqref{eq:psi_evolution}: the maximal return time is $\frac{2\pi}{\lambda_{min}}$.
Now take a matrix $\psi(0)$ in the the slice $S_0$, and let $\psi(t)$ denote the solution to Equation~\eqref{eq:psi_evolution}.

\medskip

\noindent
{\bf Claim: } Each crossing is regular and contributes positively.

\medskip

To see this, recall that the crossing form of a path $\psi$ in $Sp(2)$ at a crossing $t$ is defined as the bilinear form
\[
\omega_0(\cdot, \dot \psi(t)  \cdot)|_{\ker(\psi(t)-\id)}.
\]
Since $U$, $V$ is a symplectic frame, we have with $Z=(u,v)$ (i.e.\ $X=uU+vV$), the following inequality:
$$
\omega_0(Z,\dot \psi Z)=d\alpha(X, \dot X) \geq \lambda_{min} (u^2+v^2),
$$
by Equation~\eqref{eq:angle_change}. This establishes the  claim.

Let $t_\ell$ denote the $\ell$-th return time to $S_0$ of the path $\psi$.
The symplectic path $\psi|_{[0,t_\ell]}$ is not necessarily a loop, but we can make it into a loop by connecting $\psi(t_\ell)$ to $\psi(0)$ while staying in the slice $S_0$.
We can and will do this by adding at most one crossing, which we make regular. 
Call the extension to a loop $\tilde \psi$.  The additional crossing that we may have inserted can contribute negatively.

Now use the loop axiom for the Robbin-Salamon index. This tells us that 
$$
\mu_{RS}(\tilde \psi)=2\mu_\ell(\tilde \psi)=2 \ell.
$$
By the catenation property of the Robbin-Salamon index and positivity of all but the last (potential) crossing we have
$$
\mu_{RS}(\psi)\geq 2\ell-2.
$$
Now consider a symplectic path $\psi$ of length $T$.
We can bound the winding number as
$$
\ell \geq \left\lfloor \frac{\lambda_{min}}{2\pi} T \right\rfloor
.
$$ 
With this in mind we obtain for a Hamiltonian arc $\gamma$ of length $T$, 
$$
\mu_{RS}(\gamma;\epsilon)
\geq \frac{2\lambda_{min}}{2\pi} T-4.
$$
When this Hamiltonian arc $\gamma$ is viewed as a Reeb arc $\gamma_R$ with Reeb action $T_R$, we can rewrite this bound as follows, using
$$
T_R =\int_0^T \alpha(X_\phi) dt \leq T \cdot \max \alpha(X_\phi).
$$
We find
$$
\mu_{RS}(\gamma_R;\epsilon)
\geq \frac{\lambda_{min}}{\pi \max \alpha(X_\phi) }T_R -4.
$$
\end{proof}

\begin{remark}\label{rk:posdefonxi}
Observe that the proof actually shows the stronger claim that index-positivity holds when the Hessian of $\phi$ restricted to the contact structure is positive-definite. Note also that the latter condition is not enough for dynamical convexity.

Finally we note that the bound obtained can be sharpened since the index is necessarily positive by observing that $\psi(0)=\id$, so it is a crossing and using that each crossing of the path $\psi$ contributes positively.
\end{remark}

\section{Strongly index-definite symplectic paths}\label{app:indexposODEs}

In this appendix, we prove a crucial index growth estimate needed in order to rule out non-relevant boundary orbits via index considerations (needed in Lemma \ref{indexpositivitylemma} in the main body of the paper).

\begin{definition}\label{def:indexposODE}
Consider the linear ODE $\dot \psi(t)=A(t)\psi(t)$, where $A:\R_{\geq 0}\to \mathfrak{sp}(2n)$ and $A(0)=0$.
Its solution is a path of symplectic matrices with $\psi(0)=\mathds{1}$. We say that the ODE is \emph{strongly index-definite} if there exist constants $c>0,d \in \mathbb{R}$, such that
$$
\vert\mu_{RS}(\psi\vert_{[0,t]})\vert\geq ct+d,
$$
where $\mu_{RS}$ is the Robbin--Salamon index \cite{RS93}.
\end{definition}

Note that we make no non-degeneracy assumptions on the symplectic paths in the above definition.

\medskip

We now consider the specific family of linear ODEs $\dot \psi(t)=A(t)\psi(t)$, where the matrix $A$ has the special form

\[
A(t)=\left(\begin{array}{@{}c|c@{}}
  R(t)
  & \begin{matrix}
  \overline{X}(t) & 0 \\
  \overline{Y}(t) & 0 \\
  \end{matrix} \\
\hline
  \begin{matrix}
  0 & 0  \\
 \overline{Y}(t) & -\overline{X}(t) 
  \end{matrix} &
  \begin{matrix}
   a(t) & 0 \\
  b(t) & -a(t)
  \end{matrix}
  \end{array}\right) \in \mathfrak{sp}(2n).
\]

Here, we use the notation $(\overline{X},\overline{Y})=(X_1,Y_1,\dots,X_{n-1},Y_{n-1})$, and we assume $R(t) \in \mathfrak{sp}(2n-2)$, $A(0)=0$. 

\begin{lemma}\label{App:indexposestimate}
Assume that the linear ODE $\dot{M}(t)=R(t)M(t)$ is strongly index-definite as an ODE in dimension $2n-2$. Then the same holds for the linear ODE $\dot\psi(t)=A(t)\psi(t)$.
\end{lemma}

\begin{proof} One may check that
$$
\mathfrak{g}=\left\{ \left(\begin{array}{@{}c|c@{}}
  R
  & \begin{matrix}
  \overline{X} & 0 \\
  \overline{Y} & 0 \\
  \end{matrix} \\
\hline
  \begin{matrix}
  0 & 0  \\
 \overline{Y} & -\overline{X} 
  \end{matrix} &
  \begin{matrix}
    a & 0 \\
  b & -a
  \end{matrix}
\end{array}\right): R \in \mathfrak{sp}(2n-2)\right\}
$$
is a Lie subalgebra of $\mathfrak{sp}(2n)$. The corresponding Lie subgroup of $Sp(2n)$ is
$$
G=\left\{ \left(\begin{array}{@{}c|c@{}}
  M
  & \begin{matrix}
  \overline{x} & 0 \\
  \overline{y} & 0 \\
  \end{matrix} \\
\hline
  \begin{matrix}
  0 & 0  \\
 \overline{u} & \overline{v} 
  \end{matrix} &
  \begin{matrix}
    \alpha & 0 \\
  \beta & \alpha^{-1}
  \end{matrix}
\end{array}\right): M \in Sp(2n-2),\; \alpha>0,\; (-\overline{y},\overline{x})\cdot M +\alpha \cdot (
 \overline{u},\overline{v})=0 \right\}
$$
We deduce that $\psi \in G$. We then write
$$
\psi=\left(\begin{array}{@{}c|c@{}}
  M
  & \begin{matrix}
  \overline{x} & 0 \\
  \overline{y} & 0 \\
  \end{matrix} \\
\hline
  \begin{matrix}
  0 & 0  \\
 \overline{u} & \overline{v} 
  \end{matrix} &
  \begin{matrix}
    \alpha & 0 \\
  \beta & \alpha^{-1}
  \end{matrix}
\end{array}\right)\in G,
$$
where $M$ is a solution to $\dot M=R M$, and consider the following homotopy of paths: 
$$
\psi_s=\left(\begin{array}{@{}c|c@{}}
  M
  & \begin{matrix}
  s\overline{x} & 0 \\
  s\overline{y} & 0 \\
  \end{matrix} \\
\hline
  \begin{matrix}
  0 & 0  \\
 s\overline{u} & s\overline{v} 
  \end{matrix} &
  \begin{matrix}
    \alpha & 0 \\
  \beta & \alpha^{-1}
  \end{matrix}
\end{array}\right).
$$
Note that $\psi_s$ is a path in $G\subset Sp(2n)$ for every $s$, and $\psi_0$ has no off-diagonal terms.
For any given $t$, this gives a homotopy in $G$ relative endpoints of $\psi\vert_{[0,t]}$ to a concatenated path of the form $\psi_0\vert_{[0,t]} \# \phi_t$, where $\phi_t(s)=\psi_s(t)$. We therefore have
\begin{equation}\label{sumformula}
\mu_{RS}(\psi\vert_{[0,t]})=\mu_{RS}(\psi_0\vert_{[0,t]})+\mu_{RS}(\phi_t).
\end{equation}
On the other hand, from the block decomposition of $\psi_0$ and the fact that the lower-block can be homotoped to a symplectic shear by joining $\alpha(t)$ to $1$, we have
\begin{equation}\label{psi0formula}
\mu_{RS}(\psi_0\vert_{[0,t]})=\mu_{RS}(M\vert_{[0,t]})\pm\frac{1}{2}\mbox{sign}(\beta(t)),
\end{equation}
where the sign depends on conventions. Moreover, one may easily check that the characteristic polynomial of an element in $G$ is completely independent of the off-diagonal terms. In particular, we obtain that
$$
\det(\psi_s-\mathds{1})=\det(\psi_0-\mathds{1})=\det(M-\mathds{1})(\alpha-1)(\alpha^{-1}-1),
$$
is independent of $s$. In other words, $\psi(t)$ is an intersection point with the Maslov cycle if and only if $\psi_0(t)$ is, and the eigenvalue $1$ has the same algebraic multiplicity for both such intersections.  Moreover, if $\psi(t)$ is not an intersection, then $\phi_t$ does not intersect the Maslov cycle at all.

\medskip
One may check that if $\alpha(t)\neq 1$, then the geometric multiplicity of $1$ as an eigenvalue of $\phi_t(s)$ is independent of $s$ (and therefore $\mu_{RS}(\phi_t)=0$ for such $t$). If $\alpha(t)=1$, this may not necessarily still hold. However, we may appeal to the following general fact, whose proof was provided to the authors by Alberto Abbondandolo:

\begin{lemma}\label{univbound}
There exists a universal bound $C=C(n)$ (depending only on dimension), such that, if $\phi:[0,1] \rightarrow Sp(2n)$ is a continuous path of symplectic matrices for which the algebraic multiplicity of the eigenvalue $1$ of the matrix $\phi(t)$ is independent of $t$, then 
$$
|\mu_{RS}(\phi)|\leq C.
$$
\end{lemma}

\begin{proof}[Proof of Lemma \ref{univbound}]$\;$

\medskip

\textbf{Step 1.} We first reduce to the case where $\phi$ has $1$ as the only eigenvalue. We have a continuous symplectic splitting $\mathbb{R}^{2n}=V(t)\oplus W(t)$ where $V(t)$ is the generalized eigenspace of $\phi(t)$ corresponding to $1$, and $W(t)$ is the direct sum of the generalized eigenspaces of $\phi(t)$ corresponding to the other eigenvalues (here, the dimensions of $V(t)$ and $W(t)$ are $t$-independent by assumption), for which $\phi(t)=\phi_V(t) \oplus \phi_W(t)$ splits symplectically. Since $\phi_W$ does not intersect the Maslov cycle by construction, we have $\mu_{RS}(\phi)=\mu_{RS}(\phi_V)+\mu_{RS}(\phi_W)=\mu_{RS}(\phi_V).$

\medskip

\textbf{Step 2.} A loop $\phi$ of symplectic matrices having $1$ as the only eigenvalue is nullhomotopic in $Sp(2n)$, and hence $\mu_{RS}(\phi)=0$. This follows for instance by the interpretation of the Robbin--Salamon index as the total winding number of the Krein-positive eigenvalues on the unit circle (see e.g.\ \cite[Lemma 1.3.7]{A01}).

\medskip

\textbf{Step 3.} The identity matrix may be joined to any symplectic matrix $M$ satisfying spec$(M)=\{1\}$ via a path $M(t)$ satisfying spec$(M(t))=\{1\}$, and for which $|\mu_{RS}(M(t))|\leq C$ for some universal bound $C$. Indeed, we may write $M=e^{JS}$ where $S$ is a symmetric matrix having $0$ as the only eigenvalue, and consider the path $M(t)=e^{tJS}$. This satisfies the required properties since $M(t)$ changes strata of the Maslov cycle only at $t=0$, the geometric multiplicity of $1$ jumping from $2n$ at $t=0$ to perhaps a lower one at $t>0$, and so the contribution of this wall-crossing to $\mu_{RS}(M)$ is universally bounded.

\medskip

The proof finishes by combining the previous steps, where we join the endpoints of a path $\phi$ as in Step $1$ to the identity as in Step $3$, use the concatenation property of $\mu_{RS}$, and appeal to Step $2$. \end{proof}

Combining Equations \eqref{sumformula} and \eqref{psi0formula} with Lemma \ref{univbound}, we conclude that
$$
|\mu_{RS}(\psi\vert_{[0,t]})-\mu_{RS}(M\vert_{[0,t]})|\leq C
$$
for some universal constant $C=C(n)$, from which the conclusion of Lemma \ref{App:indexposestimate} is immediate.\end{proof}

\end{document}